# MEAN FIELD CONVERGENCE OF A MODEL OF MULTIPLE TCP CONNECTIONS THROUGH A BUFFER IMPLEMENTING RED

By D. R. McDonald[1] and J. Reynier[2]

*University of Ottawa and École Normale Supérieure*

RED (Random Early Detection) has been suggested when multiple TCP sessions are multiplexed through a bottleneck buffer. The idea is to detect congestion before the buffer overflows by dropping or marking packets with a probability that increases with the queue length. The objectives are reduced packet loss, higher throughput, reduced delay and reduced delay variation achieved through an equitable distribution of packet loss and reduced synchronization.

Baccelli, McDonald and Reynier [*Performance Evaluation* **11** (2002) 77–97] have proposed a fluid model for multiple TCP connections in the congestion avoidance regime multiplexed through a bottleneck buffer implementing RED. The window sizes of each TCP session evolve like independent dynamical systems coupled by the queue length at the buffer. The key idea in [*Performance Evaluation* **11** (2002) 77–97] is to consider the histogram of window sizes as a random measure coupled with the queue. Here we prove the conjecture made in [*Performance Evaluation* **11** (2002) 77–97] that, as the number of connections tends to infinity, this system converges to a deterministic mean-field limit comprising the window size density coupled with a deterministic queue.

**1. Introduction.** Imagine the scenario where $N$ work stations in a university department are connected by a switched ethernet to a departmental router. If every work station simultaneously FTPs a file to some distant machine, then the output buffer in the router will be a bottleneck. We study the interaction of $N$ TCP/IP connections in the congestion avoidance phase of TCP Reno routed through a bottleneck queue.

Upon receiving a TCP packet, the recipient sends back an acknowledgment packet so there is one Round Trip Time (RTT) between the time a

Received June 2003; revised January 2005.
[1]Supported in part by NSERC Grant A4551.
[2]Research initiated during an internship at the University of Ottawa.

packet is sent and the acknowledgment is received. The acknowledgment contains the sequence number of the highest value in the byte stream successfully received up to this point in time. By counting the number of packets sent but not yet acknowledged, each source implements a window flow control which limits the number of packets from this connection allowed into the network during one RTT. Duplicate acknowledgments are generated when packets arrive out of order or when a packet is lost.

The link rate of the router is $NL$ packets per second. We assume packets have equal mean sizes of 1 data unit. We assume the packets from all connections join the queue at the bottleneck buffer and we denote by $Q_N(t)$ the average queue per flow. We assume the scheduling to be FIFO.

We imagine the source writes its current window size and the current RTT in each packet it sends, where by the current RTT we mean the RTT of the last acknowledged packet:

- $W_n^N(t)$ is defined to be the window size written in a packet from connection $n$ arriving at the router at time $t$.
- $R_n^N(t)$ is defined to be the RTT written in a packet from connection $n$ arriving at the router at time $t$ (this RTT is the sum of the propagation delay plus the queueing delay in the router).
- We shall assume that connections can be divided in $d$ classes $K_c$, where $c \in [1, 2, \ldots, d]$ where the transmission time $T_n$ of any connection $n \in K_c$ equals the common transmission time $T_c$ for class $c$ (the notation will be clear from context).

At time $t$, source $n$ has sent $W_n^N(t)$ packets into the network over the last $R_n^N(t)$ seconds. The acknowledgments for these packets arrive at rate $W_n^N(t)/R_n^N(t)$ on average. New packets are being sent at the rate acknowledgments come back to the source so we will define the instantaneous transmission rate of source $n$ at time $t$ to be $X_n^N(t) = W_n^N(t)/R_n^N(t)$. This definition models the transmissions over any long period of time $T$. During time $T$, the total number of packet-minutes of work done by the network for connection $n$ is

$$\int_0^T W_n^N(t)\,dt = \int_0^T X_n^N(t) R_n^N(t)\,dt.$$

Consequently, our definition of the transmission rate is consistent with a Little type formula which calculates the work done as the integral of the packet arrival rate $X_n^N(t)$ times the work done for each packet, $R_n^N(t)$.

Under TCP Reno, established connections execute congestion avoidance where the window size of each connection increases by one packet each time a packet makes a round trip, that is, each $R_n^N$ as long as no losses or timeouts occur. During this phase, the rate the window of connection $n$ increases is approximately $1/R_n^N$ packets per second. The only thing restraining the growth of transmission rates is a loss or timeout. When a loss occurs the window is reduced by half.

The source detects a loss when three duplicate acknowledgments arrive. The source cuts the window size in half and then starts a fast retransmit/fast recovery by immediately resending the lost packet. Fast retransmit/fast recovery ends and congestion avoidance resumes when the acknowledgment of the retransmitted packet is received by the source. We assume the losses are only generated by the RED (Random Early Detection) active buffer management scheme or by tail-drop. We neglect the possibility of transmission losses.

We also neglect the possibility that some of the connections fall into timeout. This may occur if there is a loss when the window size is three or less. In this case there can't be three duplicate acknowledgments. The source can't recognize that a loss has occurred and essentially keeps on waiting for a long timeout period. Alternatively, if a retransmitted packet is lost, the source will fall into timeout. Two losses in the same RTT may not produce a timeout and may have a different effect on the window size depending on the version of TCP being used. The detection of the second loss may provoke a retransmission, but no window reduction with NewReno [11] or with SACK [13] or will only provoke a second window reduction after the acknowledgment of the retransmission of the first lost packet; that is, the discovery of the second loss occurs more than one round trip time after the second loss occurred. When the timeout period elapses, the source restarts quickly using slow-start and attempts to re-enter the congestion avoidance phase. Losses which occur simultaneously with packets arriving out of order are also a major cause of timeouts. In practice, a certain proportion of the connections will be in timeout at any given time. In effect, one has to redefine $N$ if one wants to compare theoretical predictions with simulations (see [2]).

We will assume the *large* buffer holds $B$ packets and that, once this buffer space is exhausted, arriving packets are dropped. Such tail-drops come in addition to the RED mechanism. Here we take the drop probability of RED (of an incoming packet before being processed) to be a function of the queue size which is zero for a queue length below $Q_{\min}$ but rises linearly to $p_{\max}$ at $Q_{\max}$ and is equal to 1 above $Q_{\max}$.

Note that this is not exactly as originally specified in [12], where the drop probability was taken as a linear function of the exponential moving average of the queue size. Since we will let the number of sources tend to infinity, this averaging out of fluctuations is not necessary and, in fact, is deleterious since it adds further delays into the system.

If all $N$ connections are in congestion, avoidance, we can reformulate this drop probability in terms of $Q^N$, the queue size divided by $N$, as $F(Q^N(t))$, where $F$ is a distribution function which is zero below $q_{\min} = Q_{\min}/N$ but rises linearly to $p_{\max}$ at $q_{\max} = Q_{\max}/N$ and jumps to 1 at $q_{\max}$. Of course, the tail-drop scheme can be considered as the limiting case when $q_{\min} = 0$, $q_{\max} = b$ ($B = Nb$) and $p_{\max} = 0$.

In [2] we used a fluid description of the queue and a continuous approximation of the loss rate of each connection to construct a model for the evolution of the fluid queue as a function of the histogram of window sizes. This model generalized the model in [15]. The goal of this paper is to prove the conjecture in [2] that the histogram or empirical measure $M_c^N(t, dw)$ of the window sizes in any class $c$ ($d = 1$ in [2]) converges to a deterministic mean field limit with measure $M_c(t, dw)$ at time $t$ and, moreover, the relative fluid queue size $Q^N(t)$ converges to a deterministic fluid queue $Q(t)$.

Construct a sequence $f^k$ of bounded, positive continuous convergence determining functions on $[0, \infty)$ (see pages 111 and 112 in [10]). Define a metric for weak convergence for probability measures on $[0, \infty)$ by defining the distance between probabilities $\mu$ and $\nu$ as

$$\|\mu - \nu\|_w := \sum_{k=1}^{\infty} 1 \wedge |\langle f^k, \mu\rangle - \langle f^k, \nu\rangle| \frac{1}{2^k},$$

where $\langle f, \mu\rangle = \int_0^\infty f(w)\mu(dw)$. Also, let $\|\cdot\|_s$ denote the Skorohod distance between two elements of $D[0, T]$.

THEOREM 1. *Under Assumptions 1 and 2 given in Section 2, the random measure of the window sizes of connections in each class $c$ converge in probability to a deterministic measure $M_c(t, dw)$; that is, $\|M_c^N(t, dw) - M_c(t, dw)\|_w \to 0$ in probability as $N \to \infty$. $M_c(t, dw)$ is the marginal distribution of $M_c(s - R_c(s), dv; s, dw)$, the deterministic limit of the joint distribution of the window sizes at time $t$ and at time $t - R_c(t)$.*

*Let $\mathcal{G} = \{g \in C_b^1(\mathbb{R}^+) : g(0) = 0\}$, where $C_b^1(\mathbb{R}^+)$ is the space of bounded functions with bounded derivatives. For $g \in \mathcal{G}, c = 1, \ldots, d$,*

$$\begin{aligned}
&\langle g, M_c(t)\rangle - \langle g, M_c(0)\rangle \\
&\quad = \int_0^t \biggl[\frac{1}{R_c(s)}\biggl\langle \frac{dg(w)}{dw}, M_c(s, dw)\biggr\rangle \\
&\qquad + \langle (g(w/2) - g(w))v, M_c(s - R_c(s), dv; s, dw)\rangle \\
&\qquad \times \frac{1}{R_c(s - R_c(s))} K(s - R_c(s))\biggr] ds \qquad (1.1)
\end{aligned}$$

$$\begin{aligned}
&\quad = \int_0^t \biggl[\frac{1}{R_c(s)}\biggl\langle \frac{dg(w)}{dw}, M_c(s, dw)\biggr\rangle \\
&\qquad + \langle (g(w/2) - g(w)), e(s, s - R_c(s), w) M_c(s, dw)\rangle \\
&\qquad \times \frac{1}{R_c(s - R_c(s))} K(s - R_c(s))\biggr] ds, \qquad (1.2)
\end{aligned}$$

*where* $\langle g, M_c(t)\rangle = \int_{w=0}^{\infty} g(w)M_c(t,dw)$, *where* $R_c(t) = T_c + Q(t - R_c(t))/L$, *where* $K(t) = F(Q(t))$ *for* $Q(t) < q_{\max}$ *and where*

$$e(s, s - R_c(s), w) = \left\langle v, \frac{M_c(s - R_c(s), dv; s, dw)}{M_c(t, dw)} \right\rangle$$

*is the conditional expectation of the window one RTT in the past, given the window is now* $w$.

*Moreover, the queue size converges in probability to a deterministic limit* $Q(t)$ *satisfying*

$$\begin{aligned} \frac{dQ(t)}{dt} &= \sum_{c=1}^{d} \kappa^c \langle w, M_c(t,dw)\rangle \frac{(1-K(t))}{R_c(t)} - L \\ &\quad - \left( \sum_{c=1}^{d} \kappa^c \langle w, M_c(t,dw)\rangle \frac{(1-K(t))}{R_c(t)} - L \right)^{-} \chi\{Q(t) = 0\}. \end{aligned} \tag{1.3}$$

*When* $Q(t) = q_{\max}$, $K(t)$ *is determined by*

$$K(t) = \max\left( p_{\max}, \frac{1}{L} \left( \sum_{c=1}^{d} \kappa^c \langle w, M_c(t,dw)\rangle \frac{(1-K(t))}{R_c(t)} \right)^{-1} \right),$$

*where* $\langle w, M_c(t,dw)\rangle = \int_w wM_c(t,dw)$ *is Lipschitz continuous in* $t$ *for each* $c$.

The numerical evaluation of the above equations is considered in Section 7.

Note that there may be a discontinuity in $K(t)$ when $Q(t)$ hits $q_{\max}$. The problem at $q_{\max}$ arises because $F$ is not continuous and certainly not Lipschitz at this point. To justify this definition, we consider the *Gentle RED* variant [18], where $b > 2q_{\max}$, and we extend the definition of $F$ to rise linearly from $p_{\max}$ to 1 between $q_{\max}$ to $2q_{\max}$ so $F$ is Lipschitz. In Section 5 we modify Gentle RED so that the drop probability increases linearly from $p_{\max}$ at $q_{\max}$ to 1 at $q_{\max} + \delta \le b$. The weak limit of this modified Gentle RED as $\delta \to 0$ gives the discontinuous $K(t)$ above. The reflection problem at queue size zero can be solved by the Skorohod construction.

It is important to emphasize that we have not justified that the fluid models in [15] or [2] are the limit of some discrete packet level model. We are proving far less; that is, that the intuitively attractive fluid models in [15] or [2] do converge to a mean field limit. This convergence is implicitly assumed in the engineering literature and the resulting limit processes are used to analyze the stability of various Active Queue Management (AQM) control strategies (RED among others) [15]. We also make the simplifying modeling assumptions made in [2], although more precise alternatives are suggested (see the Doppler factor, $[1 - \frac{d}{dt}R_n^N(t)]$ in (2.1)). These simplifications (or the assumption that mean windows one RTT apart are uncorrelated made in [15] but not by [2]) may be too gross and, to date, nobody has done the network measurements to check which assumptions make a significant difference. This is a major failing because the control theoretic analysis, as proposed, for instance, in [15] or [8], may be highly sensitive to these assumptions.

The conclusion, as far as RED is concerned, is negative. If the RED parameters are not chosen properly (as a function of RTT), then RED is unstable. The fixed point described in [2] may be unstable and a tiny oscillation is amplified until the queue size oscillates wildly. Even with a large buffer, the utilization may drop to less than one. Although one may criticize lacunae in the model, this conclusion is verified by simulation and for this reason, RED is rarely activated even though it is implemented on most routers.

We are mainly interested in a mathematical proof of the convergence to the mean field so we will ignore timeouts and slow-start, as well as special details of congestion avoidance which would only serve to obscure the main ideas. We will, nevertheless, sketch how these extensions could be handled. Our method could be adapted to proving mean field limits for control schemes other than RED. [19] and [7] analyze time slotted rate-based and queue-based models with delay where the number of sources tend to infinity. For their queue-based model, they prove the convergence of the queue to a deterministic limit and propagation of chaos; that is, the transmission rates of each source converge to a system driven by the deterministic queue size. In their model the limiting, deterministic queue is not coupled with limiting distribution of the window sizes and there is no associated partial differential equation as in [2]. The only comparable analysis is the recent thesis [3].

The structure of the paper is as follows. In Section 2 we model the system of $N$ windows coupled with the queue and then formulate this as a histogram of window sizes coupled with a queue. In Section 2.3 we summarize the mean field limit. The proof of the existence of this limit follows in Section 3. In Section 4 we establish the convergence to a unique limit. Finally, in Section 6 we establish Theorem 1.

The TCP model we are studying can be viewed as $N$ dynamical systems [the $N$ window sizes $\mathbf{W}^N := (W_1^N, \ldots, W_N^N)$] which evolve independently except through a shared resource (the queue $Q^N$). The dynamics of the shared resource depend only on the distribution of the dynamical systems (in this case on the average window size). The standard approach is to prove existence and then uniqueness of the limit. That's what we do here, but the mathematical innovation is to first create a modified system [here $(\mathbf{W}^N, Q^N)$ is modified to $(\mathcal{W}^N, \mathcal{Q}^M)$] where the dynamics of the shared resource (the modified queue $\mathcal{Q}^N$) depend on the expected value of the distribution (or, in this case, the expected value of the average window sizes of the modified system). Essentially, we just stick an expectation in front of the interaction term [here we replace the average window size $\overline{\mathcal{W}}^N(t)$ by $E\overline{\mathcal{W}}^N(t)$].

Since the modified shared resource $\mathcal{Q}^N$ is deterministic, the modified dynamical systems are independent. Moreover, it is easy to pick a convergent subsequence for the shared resource (here for $\mathcal{Q}^N \to \mathcal{Q}$). It is then easy to prove $\mathcal{W}_n^N$ converges to a limit $\mathcal{W}_n$ along the subsequence for each component $n$. This gives the existence of an infinite modified system [here $(\mathcal{W}, \mathcal{Q})$, where $\mathcal{W} = (\mathcal{W}_1, \mathcal{W}_2, \ldots)$]. Next, the key remark is that by the law of large numbers (and boundedness),

$$\overline{\mathcal{W}}(t) := \lim_{N\to\infty} \frac{1}{N}\sum_{n=1}^{N} \mathcal{W}_n(t) = E\overline{\mathcal{W}}(t);$$

that is, the infinite modified system is, in fact, a limit of the original system! Here this means we can rewrite $(\mathcal{W}, \mathcal{Q})$ as $(\mathbf{W}, Q)$ with $\mathbf{W} := (W_1, W_2, \ldots,)$ since the interaction term is $\overline{\mathbf{W}}(t)$; that is, the interaction is through the window average and not the expected value of the window average.

Next, we use a coupling argument to show each original dynamical system (here each $W_n^N$) converges almost surely in Skorohod norm to the infinite limit (here $W_n^N$ converges to $W_n$). This proves the propagation of chaos where each dynamical system $(W_n)$ is independent and interacts with the other systems only through the deterministic shared resource $Q$. From this, we can show the mean field convergence (here the convergence of $Q^N$ to $Q$ and the histogram of the $\mathbf{W}^N$ to the mean field limit).

We have used Kurtz's approach (e.g., [9]) of *bringing back the particles*; that is, not projecting onto the histogram of window sizes as is standard (e.g., [5]). We believe this is the most effective way to handle feedback delay because the histogram of window sizes is not a state. We think our approach has potential in other contexts.

## 2. The $N$-particle system and mean-field limit.

2.1. *The $N$-particle Markov process.* Our model takes into account the delay of one round trip time between the time the packet is killed and the time when the buffer receives the reduced rate. We assume window reductions at connection $n$ occur because of a loss one round trip time in the past. To first order, the probability of a window reduction between time $t$ and $t+h$ is

$$\begin{aligned} &\int_{t-R_n^N(t)}^{t+h-R_n^N(t+h)} \frac{W_n^N(s)}{RTT_n^N(s)} K^N(s)\,ds \\ &\quad \sim \left[1 - \frac{d}{dt}R_n^N(t)\right] \frac{W_n^N(t-R_n^N(t))}{R_n^N(t-R_n^N(t))} K^N(t-R_n^N(t))h, \end{aligned} \tag{2.1}$$

since the probability a packet is dropped is proportional to $W_n^N(t-R_n^N(t))/R_n^N(t-R_n^N(t))$, the transmission rate one RTT in the past, times $K^N(t-R_n^N(t))$, the drop probability one RTT in the past. The Doppler term $[1-\frac{d}{dt}R_n^N(t)]$ is a small correction that was overlooked in [2] and we will ignore it here.

There are many ways of actually implementing packet drops once the drop probability $p=K^N(t)$ is determined at time $t$. We could drop packets deterministically one every $1/p$ packets, but this may introduce unwanted synchronization. In fact, [12] proposes two methods. In the first we simply generate a Bernoulli random variable with probability $p$ of dropping a packet. In the second the dropped packet is chosen uniformly among the next $1/p$ packets. In fact, it won't matter which method is employed since, as $N\to\infty$, the contribution of each flow becomes negligible. Consequently, the packet arrivals of connection $n$ are enormously spread out among the other packets. As far as connection $n$ is concerned, packets are dropped randomly with probability $p=K^N(t)$ at time $t$. We therefore model the process of window reductions by a Poisson point process with stochastic intensity

$$\lambda_n^N(t) := \frac{W_n^N(t-R_n^N(t))}{R_n^N(t-R_n^N(t))} K^N(t-R_n^N(t))$$

[we can assume $W_n^N(t)=w_n$ for $t<0$].

Of course, the second method proposed in [12] would induce a weak dependence between the Poisson processes for different connections. However, the interaction between flows is via the average window size and the minor weak dependence won't prevent the average from converging to a deterministic limit. We will assume the first method is used, but it would be possible to alter the argument in Sections 3 and 6 to account for weak dependence.

We can construct the simple point process of window reductions:

$$N_n^N(t) = \int_0^t \int_0^\infty \chi_{[0,\lambda_n^N(v)]}(u) \Upsilon_n^N(du,dv),$$

where the $\Upsilon_n^N(u,v)$ are two-dimensional Poisson processes with intensities 1 on $[0,T]\times[0,\infty)$. In addition, the sources evolve independently given the trajectory of $K^N$. We therefore assume the $\Upsilon_n^N(u,v)$ are independent. In order to derive strong convergence theorems as $N\to\infty$, we shall suppose, without loss of generality, $\Upsilon_n^N(u,v)=\Upsilon_n(u,v)$, where $\{\Upsilon_n\}_{n=1,\dots,\infty}$ is a sequence of i.i.d. two-dimensional Poisson processes with intensity 1 defined on a probability space $(\Omega,\mathcal{F},P)$. In fact, Poisson processes in $\{(t,u):0\le u\le \overline{\lambda}(t), 0\le t\le T\}$ would do where $\overline{\lambda}$ is defined after Assumption 1 as an a priori bound on the transmission rate. Define $\mathcal{F}_t=\sigma\{\Upsilon_n(u,v); v\le t, n=1,2,\dots\}$.

The above is a different version of the point process of window reductions than that in [2]. The laws are the same, so the resulting dynamical systems have the same distribution. Consequently, the convergence in probability proved in Theorem 1 is also valid for the version used in [2].

*Differential equation for queue size.* For $Q(t) < q_{\max}$, $K^N(t) = F(Q^N(t))$ and the rate of change of the fluid buffer is given by

$$N\frac{dQ^N(t)}{dt} = \sum_{n=1}^{N} \frac{W_n^N(t)}{R_n^N(t)}(1 - K^N(t)) - NL + \left(\sum_{n=1}^{N} \frac{W_n^N(t)}{R_n^N(t)}(1 - K^N(t)) - NL\right)^{-} \chi\{Q^N(t) = 0\}$$

since the proportion $K^N(t) := F(Q^N(t))$ of the total fluid,

$$\sum_{n=1}^{N} \frac{W_n^N(t)}{R_n^N(t)},$$

is lost. The second term prevents the queue size from becoming negative. In effect, the queue can stick at 0 until a sufficient number of connections increase their window size.

Dividing by $N$ gives

$$\frac{dQ^N(t)}{dt} = \frac{1}{N}\sum_{n=1}^{N} \frac{W_n^N(t)}{R_n^N(t)}(1 - K^N(t)) - L + \left(\frac{1}{N}\sum_{n=1}^{N} \frac{W_n^N(t)}{R_n^N(t)}(1 - K^N(t)) - L\right)^{-} \chi\{Q^N(t) = 0\}, \tag{2.2}$$

with $Q^N(0) = q(0)$.

If $Q^N(t)$ reaches $q_{\max}$ and

$$(1 - p_{\max})\frac{1}{N}\sum_{n=1}^{N} \frac{W_n^N(t)}{R_n^N(t)} > L,$$

then the queue must jitter at $q_{\max}$ and the loss probability $K^N(t)$ is determined by

$$(1 - K^N(t))\frac{1}{N}\sum_{n=1}^{N} \frac{W_n^N(t)}{R_n^N(t)} = L.$$

In other words, if $Q^N(t) = q_{\max}$, then

$$K^N(t) = \max\left\{p_{\max}, 1 - L\left(\frac{1}{N}\sum_{n=1}^{N} \frac{W_n^N(t)}{R_n^N(t)}\right)^{-1}\right\}.$$

To justify the above definition of $K^N(t)$, one would really need to show the loss probability of a packet model jittering a $q_{\max}$ converges weakly to $K^N(t)$. Instead, we show the loss probability as *Gentle RED* converges to $K^N(t)$ as Gentle RED converges to RED (see Section 5). We should also contrast this loss rate with that for small buffers (i.e., $B$ is constant as $N \to \infty$) studied by [17]. For a small buffer, fluctuations will cause packet losses long before the total transmission rate reaches the link rate $NL$. Essentially one can model $K^N$ as $L_B(\frac{1}{LN}\sum_{n=1}^N \frac{W_n^N(t)}{R_n^N(t)})$, where $L_B$ can be calculated by finding the equilibrium distribution of a suitable Markov chain as in [16]. However, since our buffer is scaled with $N$, fluctuations like this can be ignored. It is worth noting that our method would allow us to prove mean field convergence for the small buffer case. There would be no equation for the queue and the round trip times are constant.

*Differential equation for windows.* There are three separate phases: congestion avoidance, timeout and slow start. We concentrate on describing the congestion avoidance phase. During congestion avoidance, while the queue is nonempty but less than the buffer size, the window size increases by one every time a complete window is acknowledged. In [2] and [15] this term was taken to be simply $1/R_n(t)$ (i.e., one packet increase per RTT) and we make the same approximation here. Note that this approximation ignores the fact that acknowledgments return to the sources at the link rate $NL$ when the queue is nonempty.

If the source detects the loss of a packet at time $t - R_n^N(t)$ because three duplicate acknowledgments arrive, the source cuts the current window size $W_n^N(t^-)$ by half to $W_n^N(t^-)/2$. The slow start threshold (ssthresh) is set to $H_n^N(t) = W_n^N(t^-)/2$. The source then begins fast retransmit and fast recovery. The lost packet is retransmitted and, through window inflation packets, continue to be sent as if the window size is constant [or at least the average transmission rate is consistent with a constant window size $W_n^N(t^-)/2$]. We will ignore this effect and assume the window size increases at a rate of $1/R_n(t)$ even during fast recover; that is, we don't include the term $(1 - \chi_{S_n(t)})$ in [2]. When the retransmitted packet is acknowledged, congestion avoidance resumes. Hence, the evolution of the window size in the congestion avoidance phase is described by the following stochastic differential equation:

$$dW_n^N(t) = \frac{1}{R_n^N(t)}\,dt - \frac{W_n^N(t^-)}{2}\,dN_n^N(t), \tag{2.3}$$

with $W_n^N(0) = w_n$, $n = 1, \ldots, N$, specified. Denote the vector of window sizes by $\mathbf{W}^N(t)$.

If we wished to model timeouts, we could define a function $U(W_n^N(t^-))$ equal to one if the connection falls into timeout during fast recovery and zero if not. Hence, the point process of falling into timeout is given by $U(W_n^N(t^-))\,dN_n^N(t)$. During the timeout phase, the window size is zero.

The connection is described by ssthresh and the remaining time in timeout. After the timeout phase elapses, the source enters slow-start and doubles its window size starting from one every RTT until the window size reaches ssthresh, at which time congestion avoidance restarts. If another loss is detected before reaching the congestion avoidance phase, the connection will go into timeout. During the slow-start phase, the connection is described by the window size and ssthresh. At any time $t$, a certain proportion of the $N$ connections will be in each phase. In the mean-field limit these proportions will converge to deterministic fractions. We will not show this here. In fact, we will simply ignore all the special details of fast recovery and timeouts.

*Assumptions on the initial state.*

ASSUMPTION 1. (i) Prior to time zero, the window size of connection $n$ in the $N$ connection system is a constant $w_n$, where $0 \le w_n \le W_{\max}$.

(ii) The transmission time $T_n$ of connection $n$ satisfies $T_{\min} \le T_n \le T_{\max}$ for all $n$.

(iii) $Q^N(0) = q(0)$ a constant.

*Bound $a(t)$ for the window size at time $t$.* From Assumption 1(i),

$$a(t) := W_{\max} + \frac{t}{T_{\min}} \ge w_n + \frac{t}{T_n} \ge W_n^N(t)$$

at every time $t$. The stochastic intensity for the Poisson point process of losses of connection $n$ is $\lambda_n^N(s) \le a(s)/T_{\min} =: \overline{\lambda}(t)$ for all $0 \le t \le T$.

Note that $(1 - K^N(t)) \ge 1 - p_{\max}$ as long as $Q^N(t) < q_{\max}$. If $Q^N(t) = q_{\max}$, then $(1 - K^N(t)) \ge LT_{\min}/a(t)$. Either way we have $(1 - K^N(t)) \ge (1 - k_{\max}) > 0$.

*Relation between RTT and queue size.* Define $\phi_n^N(s)$ to be the future round trip time written into a packet leaving the source at time $s$. For the above scenario, $\phi_n^N(s) = T_n + Q^N(s)/L$. Also note that $s + \phi_n^N(s)$ is monotonic because the derivative, if $Q^N(t) < q_{\max}$, is

$$\begin{aligned} 1 + \frac{1}{L}\frac{dQ^N(s)}{ds} &= 1 + \frac{1}{L}\left(\frac{1}{N}\sum_{n=1}^{N}\frac{W_n^N(t)}{R_n^N(t)}(1 - K^N(t)) - L\right) \\ &\ge \frac{1}{L}\frac{1}{N}\sum_{n=1}^{N} W_n^N(t)\frac{1}{T_{\min}}(1 - k_{\max}). \end{aligned}$$

This is positive unless all the window sizes are zero and this has probability zero. If $Q^N(t) = q_{\max}$, then the derivative is one.

Now define the RTT of connection $n$ as marked in packets arriving at the router at time $t$ by $R_n^N(t) = t - s = \phi_n^N(s)$ if $s + \phi_n^N(s) = t$. Since $s +$ $\phi_n^N(s)$ is monotonic, $R_n^N(t)$ is well defined and $\phi_n^N(t-R_n^N(t)) = R_n^N(t)$. Also, substituting $R_n^N(t) = t - s$ into $\phi_n^N(s) = T_n + Q^N(s)/L$, we get that $R_n^N$ satisfies

$$R_n^N(t) = T_n + Q^N(t - R_n^N(t))/L. \tag{2.4}$$

Moreover, by taking the derivative of (2.4), we get

$$(1 - \dot{R}_n^N(t)) = \frac{1}{1 + \dot{Q}^N(t - R_n^N(t))/L} \tag{2.5}$$

$$= L\left(\frac{1}{N}\sum_{n=1}^{N}\frac{W_n^N(t-R_n^N(t))}{R_n^N(t-R_n^N(t))}(1 - K^N(t - R_n^N(t)))\right)^{-1}. \tag{2.6}$$

2.2. *Reformulation in terms of a measure-valued process.*

*Classes of connections.* We will assume there are $d$ classes of connections $K_c$, $c = 1, \dots, d$, and all connections in class $c$ have the same transmission time $T_c$. Hence, $R_n^N = R_c^N$ for all $n \in K_c$. We will also assume the proportion of the $N$ connections in class $c$ is $\kappa_c^N$. In addition to Assumptions 1, we assume the following:

Assumption 2. Assumptions on connection classes:

(iv) The proportion of users in the class $c : \kappa_c^N \to \kappa_c$ for $c = 1, \dots, d$ as $N \to \infty$.

(v) Let $\mu_c^N$ be the histogram of windows of connections from class $c$ at time 0. We suppose that, for all $c$, $\mu_c^N$ converges weakly to $\mu_c$ as $N \to \infty$, where the support of $\mu_c$ is concentrated on $[0, W_{\max}]$.

*Measure-valued process.* In order to study the limiting behavior of the system as the number of connections $N$ goes to infinity, we will first define the empirical process (see [5, 6]) of those connections in class $K_N^c$. For any Borel set $A$, define

$$M_c^N(t, A) := \frac{1}{\kappa_c^N N}\sum_{n=1}^{N}\chi_A(W_n^N(t))\chi\{n \in K_c\} \tag{2.7}$$

to be the associated probability-measure-valued process taking values in $M_1(\mathbb{R}^+)$, the set of probability measures on $\mathbb{R}^+ = [0, \infty)$ furnished with the topology of weak convergence.

*How the future is determined.* The sequence $\{\Upsilon_n(u,v)\}_{n\in\mathbb{N}}$ of independent Poisson processes with intensity 1 was defined on a probability space $\{\Omega,\mathcal{F},P\}$. $\mathbf{W}^N(t)\equiv(W_1^N(t),\ldots,W_N^N(t))$, $Q^N(t)$ and $M^N(t)\equiv(M_1^N(t),\ldots,M_d^N(t))$ can be constructed path by path as processes defined on $\{\Omega,\mathcal{F},P\}$ taking values in $(\mathbb{R}^+)^\infty$, $\mathbb{R}^+$ and $M_1(\mathbb{R}^+)^d$, where the coordinates of $\mathbf{W}^N(t)$ above $N$ are zero. It suffices to assume $W_n^N(t)=w_n$ for $t\le 0$ and build the solution of the system (2.3), (2.2) pathwise from jump point to jump point of $\Upsilon_n$.

*Reformulating* (2.3). Let

$$\langle g,\mu\rangle=\int g(w)\mu(dw)$$

and

$$\langle Id,\mu\rangle=\int w\mu(dw),$$

so

$$\overline{W}_c^N(s):=\langle Id,M_c^N(s)\rangle=\frac{1}{\kappa_c^N N}\sum_{n=1}^N W_n^N(s)\chi\{n\in K_c\}.$$

If $g\in\mathcal{G}$, then

$$\begin{aligned}&\langle g,M_c^N(t)\rangle-\langle g,M_c^N(0)\rangle\\ &\quad=\frac{1}{\kappa_c^N N}\sum_{n=1}^N\chi\{n\in K_c\}\int_0^t\left[\frac{dg}{dw}(W_n^N(s))\frac{ds}{R_c^N(s)}\right.\\ &\qquad\qquad\left.+\left(g(W_n^N(s^-)/2)-g(W_n^N(s^-))\right)dN_n^N(s)\right].\end{aligned}\tag{2.8}$$

In Section 6 we consider the limit of the above as $N$ goes to infinity to obtain an equation for the evolution of the distribution of the windows.

*Reformulating* (2.2). For $Q^N(t)<q_{\max}$, $K^N(t)=F(Q^N(s))$ and

$$\begin{aligned}&Q^N(t)-Q(0)\\ &\quad=\int_0^t\left[\sum_{c=1}^d\kappa_c^N\langle Id,M_c^N(s)\rangle\frac{(1-K^N(s))}{R_c^N(s)}-L\right.\\ &\qquad\left.+\left(\sum_{c=1}^d\kappa_c^N\langle Id,M_c^N(s)\rangle\frac{(1-K^N(s))}{R_c^N(s)}-L\right)^-\chi\{Q^N(s)=0\}\right]ds,\end{aligned}\tag{2.9}$$

where

$$R_c^N(t)=T_c+Q^N(t-R_c^N(t))/L.\tag{2.10}$$

If $Q^N$ jitters at $q_{\max}$, then the loss probability $K^N(t)$ is given by

$$K^N(t) = \max\left\{p_{\max}, 1 - L\left(\sum_{c=1}^{d} \kappa_c^N \langle Id, M_c^N(t)\rangle \frac{(1-K^N(t))}{R_c^N(t)}\right)^{-1}\right\}.$$

2.3. *Summary of the mean-field limit.* We wish to show that $(\mathbf{W}^N(t), Q^N(t))$, the unique solution to the $N$-particle system, converges as $N\to\infty$. In Section 3 we first prove the existence of the following limit.

THEOREM 2. *If Assumptions 1 and 2 hold, then there exists a unique strong solution $(\mathbf{W}, Q, (M_1,\dots,M_d))$ to the following system. For $Q(t) < q_{\max}$, $K(t) = F(Q(t))$ and*

$$\begin{aligned} &Q(t) - Q(0) \\ &\quad = \int_0^t \Bigg[\sum_{c=1}^{d} \kappa^c \langle Id, M_c(s)\rangle \frac{(1-K(s))}{R_c(s)} - L \\ &\qquad + \left(\sum_{c=1}^{d} \kappa^c \langle Id, M_c(s)\rangle \frac{(1-K(s))}{R_c(s)} - L\right)^{-} \chi\{Q(s)=0\}\Bigg]\,ds. \end{aligned} \tag{2.11}$$

*When $Q(t) = q_{\max}$, then $K(t)$ satisfies*

$$K(t) = \max\left\{p_{\max}, 1 - L\left(\sum_{c=1}^{d} \kappa^c \langle Id, M_c(s)\rangle \frac{(1-K(t))}{R_c(t)}\right)^{-1}\right\}.$$

*Each window evolves according to*

$$dW_n(t) = \frac{1}{R_n(t)}\,dt - \frac{W_n(t^-)}{2}\,dN_n(t), \tag{2.12}$$

*where $W_n(0) = w_n$, $n = 1,\dots,$ are specified, where*

$$N_n(t) = \int_0^t \int_0^\infty \chi_{[0,\lambda_n(t)]}(v)\,d\Upsilon_n(u,v),$$

*where*

$$\lambda_n(s) = \frac{W_n(s - R_n(s))}{R_n(s - R_n(s))} K(s - R_n(s))$$

*and where $M_c(t)$ is defined by*

$$\langle g, M_c(t)\rangle = \lim_{N\to\infty} \frac{1}{\kappa_c^N N} \sum_{n=1}^{N} g(W_n(t))\chi\{n \in K_c\}$$

*and, as a consequence, we can define*

$$\overline{W}_c(s) = \langle Id, M_c(t)\rangle = \lim_{N\to\infty} \frac{1}{N\kappa_c^N} \sum_{n=1}^{N} W_n(s)\chi\{n \in K_c\}.$$

*Consequently, from* (2.11),

$$\begin{aligned} Q(t) - Q(0) = \int_0^t \Bigg[ & \sum_{c=1}^{d} \kappa_c \frac{\overline{W}_c(s)}{R_c(s)}(1 - K(s)) - L \\ & + \left( \sum_{c=1}^{d} \kappa^c \frac{\overline{W}_c(s)}{R_c(s)}(1 - K(s)) - L \right)^{-} \chi\{Q(s) = 0\} \Bigg] ds. \end{aligned} \tag{2.13}$$

*The solutions* $Q(t)$ *and* $\mathbf{R}(t) = (R_1(t), R_2(t), \dots)$ *are deterministic, as are the* $M_c(t)$*. Finally, the components of* $\mathbf{W}$ *are independent processes.*

Define

$$\overline{\mathbf{S}}_N^N(s) = \frac{1}{N} \sum_{n=1}^{N} \frac{W_n^N(s)}{R_n^N(s)} = \sum_{c=1}^{d} \kappa_c^N \frac{\overline{W}_c^N(s)}{R_c^N(s)}, \tag{2.14}$$

where $\overline{W}_c^N(s)$ is the average window size of connections in $K_c$ among the first $W_1^N, \dots, W_N^N$. Define $\overline{\mathbf{S}}_N(s)$ analogously from $\mathbf{W}$. From Theorem 2, we can define

$$\overline{\mathbf{S}}(s) = \lim_{N\to\infty} \overline{\mathbf{S}}_N(s) = \lim_{N\to\infty} \frac{1}{N} \sum_{n=1}^{N} \frac{W_n(s)}{R_n(s)}. \tag{2.15}$$

In Section 4 we will prove Theorem 3 and show that there is only one strong solution to (2.13), (2.12) and that, in fact, the solution to (2.3) and (2.2) converges to this strong solution.

THEOREM 3. *If Assumptions 1 and 2 hold, then* $\|M_c^N(t) - M_c(t)\|_w$, $\|Q^N(t) - Q(t)\|_s$ *and* $\|K^N(t) - K(t)\|_s$ *converge to zero in probability, where* $M_c(t)$*,* $q(t)$ *and* $K(t)$ *are deterministic functions of* $t \in \mathbb{R}^+$ *into* $M_1(\mathbb{R}^+)$*,* $\mathbb{R}^+$ *and* $\mathbb{R}^+$*, respectively, given in Theorem 2.*

Let $P^N$ be the measure induced on $D^d[0,T] \times C[0,T]$ by $((\overline{W}_1^N, \dots, \overline{W}_d^N), Q^N)$.

LEMMA 1. *Under Assumptions 1 and 2, the measures* $P^N$ *are tight.*

PROOF. We check the conditions for Theorem 12.3 in [4]. Condition (i) is immediate since the sequences $((\overline{W}_1^N, \dots, \overline{W}_d^N), Q^N(t))$ are bounded.

In condition (ii) we are given positive constants $\varepsilon$ and $\eta$. We can pick $\Delta$ sufficiently small that the maximum growth of a window over a duration of length $\Delta$ is less than $\varepsilon/2$; that is, pick $\Delta < T_{\min}\varepsilon/2$. Also pick $\Delta$ sufficiently small that the probability of the event $B$ that a Poisson process with intensity $\overline{\lambda}(T)$ jumps twice within a duration $\Delta$ is less than $\eta\varepsilon a(T)/2$.

Note that, by the construction of the window $W_n^N$, the event that the window is cut in half two times within a duration $\Delta$ up to time $T$ is contained in $B_n$, the event where $\Upsilon_n(t, \overline{\lambda}(T))$ has two jumps in an interval of duration $\Delta$. Also, note the worst oscillation a window can make over a duration $\Delta$ is $a(T)$; that is, the biggest drop possible. Hence, the modulus $w'_c(\Delta)$ of any trajectory of $\overline{W}_c^N, c = 1, \dots, d$, as defined at (12.6) in [4] satisfies

$$w'_c(\Delta) \leq \frac{1}{N} \sum_{n=1}^{N} (a(T)\chi_{B_n} + \varepsilon/2).$$

Hence,

$$\begin{aligned} P^N(w'_c(\Delta) \geq \varepsilon) &\leq P\left(a(T)\frac{1}{N}\sum_{n=1}^{N}\chi_{B_n} \geq \varepsilon/2\right) \\ &\leq \frac{2}{\varepsilon a(T)} P(B) \leq \eta. \end{aligned}$$

Since we can make the above estimates simultaneously for each of the $d$ classes, we see condition (ii) holds for the oscillations of $(\overline{W}_1^N, \dots, \overline{W}_d^N)$. The oscillations of $Q^N(t)$ and of $(R_1^N(t), \dots, R_d^N(t))$ are uniformly bounded in $C[0,T]$ because each trajectory is uniformly Lipschitz. This shows the $P^N$ are tight. $\square$

Using Lemma 1, we can extract a subsequence $N_k$ such that $Q^N$ and $(\overline{W}_1^N, \dots, \overline{W}_d^N)$ converge almost surely to $Q^\infty$ and $(\overline{W}_1^\infty, \dots, \overline{W}_d^\infty)$ in Skorohod norm. The convergence of the components $W_n^N$ follows. Unfortunately, we wouldn't even know that the limits $Q^\infty$ and $(\overline{W}_1^\infty, \dots, \overline{W}_d^\infty)$ are deterministic.

Using the above method and Jakubowski's criterion (cf. [6], Theorem 3.6.4), we might be able to check that the measure valued processes $M_1^N(t), \dots, M_d^N(t)$ in $D([0,T], M_1(\mathbb{R}^+))$ are tight. We might even check that the bivariate measure valued processes $M_1^N(t - R^N(t), dv; t, dw), \dots, M_d^N(t - R^N(t), dv; t, dw)$ in $D([0,T], M_1(\mathbb{R}^+)^2)$ are tight. We could then pick a convergent subsequence and carry out the analysis in Section 6. We would obtain a limiting solution to Theorem 1. Unfortunately, we won't know the solution is unique (and deterministic) because (1.1) and (1.3) don't determine the solution since $(M_1(t), \dots, M_d(t))$ and $Q(t)$ is not a state (because of the delay in the system).

It may be possible to rectify this by defining the state at time $t$ to be the entire trajectory of the measures $M_1^N(t), \ldots, M_d^N(t)$ and $Q(t)$ back at least one RTT. Indeed, the numerical procedure proposed in Section 7 shows how to maintain all this information to numerically solve (1.1) and (1.3). However, instead of trying to characterize tightness of such an ugly space, we will proceed in a more direct manner in the next section.

**3. Existence of a limit.** In this section we show the existence of the solution to (2.13), (2.12).

3.1. *Modified system.* We now introduce the modified system discussed in the introduction where $Q^N$ is forced to be deterministic by modifying the equation for the evolution of the queue to (3.1). Then we extract a deterministic limit that turns to be a limit of our initial system.

Let $\mathcal{K}^N(t) = F(\mathcal{Q}^N(t))$ if $\mathcal{Q}(t) < q_{\max}$, where $\mathcal{Q}^N$ is defined by

$$\begin{aligned} &\mathcal{Q}^N(t) - Q(0) \\ (3.1)\quad &= \int_0^t \Bigg[ \sum_{c=1}^d \kappa_c^N (E\langle Id, \mathcal{M}_c^N(s)\rangle) \frac{(1-\mathcal{K}^N(s))}{\mathcal{R}_c^N(s)} - L \\ &\quad + \left( \sum_{c=1}^d \kappa_c^N (E\langle Id, \mathcal{M}_c^N(s)\rangle) \frac{(1-\mathcal{K}^N(s))}{\mathcal{R}_c^N(s)} - L \right)^- \chi\{\mathcal{Q}(s) = 0\} \Bigg] \, ds, \end{aligned}$$

where $\mathcal{W}^N = (\mathcal{W}_1^N, \ldots, \mathcal{W}_N^N)$ satisfies the analogue of (2.3), where $\mathcal{W}_n(0) = w_n$, where

$$\mathcal{R}_c^N(t) = T_c + \mathcal{Q}^N(t - \mathcal{R}_c^N(t))/L, \tag{3.2}$$

and where

$$\langle Id, \mathcal{M}_c^N(s)\rangle = \frac{1}{\kappa_c^N N} \sum_{n=1}^N \mathcal{W}_n^N(s) \chi\{n \in K_c\}.$$

If $\mathcal{Q}^N(t) = q_{\max}$, then the loss probability $\mathcal{K}^N(t)$ satisfies

$$\mathcal{K}^N(t) = \max\left\{ p_{\max}, 1 - L \left( \sum_{c=1}^d \kappa_c^N \frac{E\langle Id, \mathcal{M}_c^N(t)\rangle}{\mathcal{R}_c^N(t)} \right)^{-1} \right\}. \tag{3.3}$$

We remark that when $\mathcal{Q}^N$ hits $q_{\max}$, the loss probability $\mathcal{K}^N(t)$ jumps from $p_{\max}$ to a value which keeps $\dot{\mathcal{Q}}^N(t) = 0$ as long as

$$\sum_{c=1}^d \kappa_c^N \frac{E\langle Id, \mathcal{M}_c^N(t)\rangle}{\mathcal{R}_c^N(t)} (1 - p_{\max}) \ge L.$$

We will show below that $E\langle Id, \mathcal{M}_c^N(t)\rangle$ is Lipschitz continuous so, in fact, $\mathcal{K}^N(t) = p_{\max}$ just before $\mathcal{Q}^N(t)$ leaves the boundary.

*Solution for a given N.* Let $t^c(0)=0$, $t^c(k+1)=t^c(k)+T_c+\mathcal{Q}^N(t^c(k))/L$ such that $t^c(k+1)-\mathcal{R}_c^N(t^c(k+1))=t^c(k)$. As long as we can define it, the sequence $(t^c(k))$ is increasing in $k$. Define $\Phi^c(t)=$ the first $k$ such that $t^c(k)>t$. We will construct our solution by recurrence from time $t_i$ to $t_{i+1}$ by defining $t_{i+1}=\min_c t^c(\Phi^c(t_i))$ starting from time $t_0=0$.

At time $t_0$, we suppose $\mathcal{W}^N(t)$ and $\mathcal{Q}^N(t)$ are given (perhaps constant) for $t\le 0=t_0$. We suppose $(\mathcal{W}^N(t),\mathcal{Q}^N(t))$ is defined for $t\le t_i$, where $t_i$ is a time such that $t_i=t^c(k)$ for some $c$ and some $k$. This is certainly true at time $t_0$. Then $\Phi^c(t_i)$ and $t^c(\Phi^c(t_i))$ are defined for all classes as is $t_{i+1}$.

Then if $t\le t_{i+1}$,

$$\lambda_n^N(t)=\frac{\mathcal{W}_n^N(t-\mathcal{R}_n^N(t))}{\mathcal{R}_n^N(t-\mathcal{R}_n^N(t))}\mathcal{K}^N(t-\mathcal{R}_n^N(t))$$

can be defined, because, for each class, $t-\mathcal{R}_c^N(t)\le t_i$ for $s\le t_{i+1}$ by the definition of $t_{i+1}$ [recall (2.6)]. Hence, the point processes $\mathcal{N}_n^N(t)$ and the trajectories $\mathcal{W}^N(t)$ are defined on $[0,t_{i+1}]$. They are bounded and measurable, thus, the expectations can be defined and, hence $\mathcal{Q}^N(t)$ can be defined. We have therefore checked the induction hypothesis up to time $t_{i+1}$.

To conclude, we need to show that $t_i\to\infty$. Notice that $t_{i+d+1}\ge\min\{t^c\times(\Phi^c(t_i)+1):c=1,\dots,d\}$ because otherwise the $d+1$ values $t_{i+j}, j=1,\dots,d+1$ must be chosen among the $d$ values $\{t^c(\Phi^c(t_i))\}$ and this is impossible. We conclude $t_{i+d+1}\ge t_i+T_{\min}$ and, therefore, $t_i\to\infty$ as $i\to\infty$.

*Lipschitz continuity of $E\langle Id,\mathcal{M}_c^N(t)\rangle$ and $\mathcal{Q}^N$ on $[0,T]$.*

$$\begin{aligned}
&|E\langle Id,\mathcal{M}_c^N(t+h)\rangle-E\langle Id,\mathcal{M}_c^N(t)\rangle|\\
&\quad\le\frac{1}{\kappa_c^N N}\sum_{n=1}^N E|\mathcal{W}_n^N(t+h)-\mathcal{W}_n^N(t)|\chi\{n\in K_c\}\\
&\quad\le\frac{h}{T_{\min}}\frac{1}{\kappa_c^N N}\sum_{n=1}^N P(\mathcal{W}_n^N(s^-)=\mathcal{W}_n^N(s),t\le s\le t+h)\chi\{n\in K_c\} \qquad (3.4)\\
&\qquad+a(T)\frac{1}{\kappa_c^N N}\sum_{n=1}^N P(\mathcal{W}_n^N(s^-)\ne\mathcal{W}_n^N(s),\\
&\qquad\qquad\text{for some } t\le s\le t+h)\chi\{n\in K_c\}. \qquad (3.5)
\end{aligned}$$

The second term arises because even multiple jumps will create a difference less than the maximum window size.

(3.4) is less than $h/T_{\min}$ and tends to zero as $h\to 0$. (3.5) is bounded by the probability a window makes a jump in an interval of length $h$. The intensity function $\lambda^N(t)$ is bounded by $\overline{\lambda}(t)=a(t)/T_{\min}$, so the probability of a jump in an interval of length $h$ is bounded by $1-\exp(a(T)h/T_{\min})$.

Hence, $E\langle Id, \mathcal{M}_c^N(t)\rangle$ is Lipschitz uniformly for $N \in \mathbb{N}$ and $t \in [0, T]$. From (3.1), it immediately follows that $\mathcal{Q}^N$ is also Lipschitz.

Let $m_c^N(t) = E\langle Id, \mathcal{M}_c^N(t)\rangle$.

*Lower bound on* $\dot{m}_c^N(t)$.

LEMMA 2. *The derivatives of* $m_c^N(t)$ *and* $\mathcal{Q}^N(t)$ *are bounded below uniformly in* $N$.

PROOF. Taking expectations,

$$
\begin{aligned}
&E\mathcal{W}_n^N(t) - w_n \\
&\qquad = E\int_0^t \left[\frac{1}{\mathcal{R}_n^N(s)}\,ds - \frac{\mathcal{W}_n^N(s^-)}{2}\,dN_n(s)\right] \\
&\qquad = \int_0^t \frac{1}{\mathcal{R}_n^N(s)}\,ds - E\int_0^t \left[\frac{\mathcal{W}_n^N(s^-)}{2}\frac{\mathcal{W}_n^N(s-\mathcal{R}_n^N(s))}{\mathcal{R}_n^N(s-\mathcal{R}_n^N(s))}\mathcal{K}^N(s-\mathcal{R}_n^N(s))\right] ds.
\end{aligned}
$$

Hence,

$$
\begin{aligned}
&\frac{dE\mathcal{W}_n^N(t)}{dt} \\
&\qquad = \frac{1}{\mathcal{R}_n^N(t)} - E[\mathcal{W}_n^N(t^-)\mathcal{W}_n^N(t-\mathcal{R}_n^N(t))]\frac{1}{2\mathcal{R}_n^N(t-\mathcal{R}_n^N(t))}\mathcal{K}^N(t-\mathcal{R}_n^N(t)) \\
&\qquad \geq \frac{1}{T_{\max}+q_{\max}/L} - \frac{E\mathcal{W}_n^N(t^-)}{2}\frac{a(t)}{T_{\min}}.
\end{aligned}
$$

Since $E\mathcal{W}_n^N(t) \leq a(t)$, it follows that $dE\mathcal{W}_n^N(t)/dt$ is strictly bounded below by $-C$, where $C$ is a positive constant that doesn't depend on $N$ or $n$. Taking the average of the $\kappa_c N$ windows with RTT, $T_c$ shows the $m_c^N(t)$ are bounded below uniformly in $N$.

Since $\mathcal{R}_n^N(t) \leq T_{\max} + q_{\max}/L$, it follows that the derivative of $\frac{1}{N}\sum_{n=1}^N \frac{EW_n^N(t)}{\mathcal{R}_n^N(t)}$ is bounded below by the same constant. The fact that $\mathcal{Q}^N(t)$ is bounded below uniformly in $N$ follows from (3.1). $\square$

LEMMA 3. *The sequence of functions* $\mathcal{K}^N(t)$ *on* $[0, T]$ *is sequentially compact.*

PROOF. As long as $\mathcal{Q}^N(t) < q_{\max}$, $\mathcal{K}^N(t) = F(\mathcal{Q}^N(t))$, so $\mathcal{K}^N(t)$ is Lipschitz uniformly in $N$. Similarly, as long as $\mathcal{Q}^N(t) = q_{\max}$, $\mathcal{K}^N(t) = 1 - L(\sum_{c=1}^d \kappa_c^N m_c^N(t)/\mathcal{R}_c^N(t))^{-1}$, and this is Lipschitz uniformly in $N$. The problem is the jumps when $\mathcal{Q}^N(t)$ hits the boundary.

We check the conditions in Theorem 12.3 in [4]. Condition (12.25) is trivial since $\mathcal{K}^N(t)$ is uniformly bounded. To check (12.26), we must show the oscillations over small intervals after excluding big jumps is as small as we like. For any $\varepsilon > 0$, we can define the set of times $J^N = \{j_i\}$ associated with jumps bigger than $\varepsilon/2$. The number of points in this set is bounded uniformly in $N$ and the spacing between the points of $J^N$ is bounded below uniformly in $N$ by some $\delta_0$. This follows immediately from Lemma 2 because, after a jump of $\mathcal{K}(t)$ when $\mathcal{Q}(t)$ hits $q_{\max}$ of size greater than $\varepsilon$, the time to decrease to $p_{\max}$ is strictly bounded below. Select a $\delta < \delta_0$ such that $\mathcal{Q}^N(t)$ and $1 - L(\sum_{c=1}^d \kappa_c^N m_c^N(t)/\mathcal{R}_c^N(t))^{-1}$ oscillate less than $\varepsilon/4$ on intervals of length less than $\delta$. Next, consider a partition $A^N$ by points $0 = t_0 < t_1 < \cdots < t_v = T$ which are $\delta$-sparse; that is, such that $\min_{1\le i\le v}(t_i - t_{i-1}) > \delta$, which includes the points in $J^N$. This is possible because these points are spaced out by more than $\delta_0$.

Now consider the maximum oscillations over any interval $[t_{i-1}, t_i)$. There are no jumps of size greater than $\varepsilon/2$ since $\mathcal{K}^N$ is right continuous and the big jumps are among the left endpoints $t_{i-1}$. Since the jumps only go up from $p_{\max}$, they don't add, so, in fact, the greatest possible oscillation is $\varepsilon/2$ for one jump plus $\varepsilon/4 + \varepsilon/4$ for the oscillations of $F(\mathcal{Q}^N(t))$ and $1 - L(\sum_{c=1}^d \kappa_c^N m_c^N(t)/\mathcal{R}_c^N(t))^{-1}$. We conclude $w'_{\mathcal{K}^N}(\delta)$, as defined in [4], is less than $\varepsilon$ and this establishes condition (12.26). $\square$

3.2. *Existence of a limit for the modified system.* In this section we shall be extracting subsequences of sequences, but we won't reflect this in our notation until the end of this subsection.

*Extraction of a limit for $\mathcal{Q}^N$ and $E\langle Id, \mathcal{M}_c^N(t)\rangle$.* $\mathcal{Q}^N(t)$ is deterministic. Moreover, the integrand in (3.1) is bounded by a constant $B$ because the window sizes up to time $t$ are bounded by $a(T)$ and RTT in greater than $T_{\min} > 0$. Hence, $\mathcal{Q}^N$ is Lipschitz uniformly for $N \in \mathbb{N}$ and $t \in [0, T]$. It follows that there is a subsequence and a Lipschitz function $\mathcal{Q}(t)$ such that $\mathcal{Q}^N(t) \to \mathcal{Q}(t)$ uniformly using the Ascoli–Arzelà theorem plus the fact that a uniform limit of Lipschitz functions is Lipschitz.

We showed above that $E\langle Id, \mathcal{M}_c^N(t)\rangle$ is Lipschitz so again using the Ascoli–Arzelà theorem we can take a further subsequence of $N$ such that, for all $c$, $m_c^N(t) = E\langle Id, \mathcal{M}_c^N(t)\rangle$ converges uniformly to a Lipschitz function $m_c(t)$.

Note that taking the limit as $N \to \infty$ in Lemma 2 gives that the derivatives of $m_c(t)$ and $\mathcal{Q}(t)$ are bounded below.

*Convergence of the RTT.* As a direct consequence of the convergence of $\mathcal{Q}^N$, $\mathcal{R}_n^N(t)$ converges uniformly to $\mathcal{R}_n(t)$, where $\mathcal{R}_n(t) = T_n + \mathcal{Q}(t - \mathcal{R}_n(t))/L$.

*Extraction of the limit for $\mathcal{K}^N(t)$.* By Lemma 3, we can extract a further subsequence so that $\mathcal{K}^N$ converges in the Skorohod topology to a limit $\mathcal{K}$ in $D[0,T]$. If $\mathcal{Q}(t) < q_{\max}$, then $\mathcal{Q}^N(t) < q_{\max}$ for $N$ large enough. Since $\mathcal{K}^N(t) = F(\mathcal{Q}^N(t))$, it follows that $\mathcal{K}(t) = F(\mathcal{Q}(t))$.

If $\mathcal{Q}(t) = q_{\max}$, we want to show $\mathcal{K}(t)$ is given by

$$\max\left\{p_{\max}, 1 - L\left(\sum_{c=1}^{d} \kappa_c \frac{m_c(t)}{\mathcal{R}_c(t)}\right)^{-1}\right\}. \tag{3.6}$$

First suppose $t$ is not a point where $\mathcal{K}$ jumps. Then there exists a small interval $I = (t - \delta_0, t]$, where $\mathcal{Q}(s) = q_{\max}$ for $s \in I$. There are two possibilities:

$$\sum_{c=1}^{d} \kappa_c m_c(t) \frac{1 - p_{\max}}{\mathcal{R}_c(t)} > L \text{ or } \sum_{c=1}^{d} \kappa_c m_c(t) \frac{1 - p_{\max}}{\mathcal{R}_c(t)} = L.$$

In the first case $\sum_{c=1}^{d} \kappa_c m_c(t) \frac{1 - F(\mathcal{Q}(t))}{\mathcal{R}_c(t)} > L$ for $t \in [t - \delta, t]$ for some $\delta$, $0 < \delta < \delta_0$, sufficiently small. Since $\mathcal{Q}^N(t)$ converges to $\mathcal{Q}(t)$, we can pick $N$ sufficiently large that $\mathcal{Q}^N(t)$ is in a narrow tube around $\mathcal{Q}(t)$. Moreover,

$$\frac{\sum_{c=1}^{d} \kappa_c^N m_c^N(t)}{\mathcal{R}_c^N(t)} (1 - F(\mathcal{Q}^N(t))) \to \frac{\sum_{c=1}^{d} \kappa_c m_c(t)}{\mathcal{R}_c(t)} (1 - F(\mathcal{Q}(t))).$$

The right-hand side is strictly greater than $L$ for $t \in [\rho - \delta, \rho]$ for $\delta$ sufficiently small and the same is true of the left-hand side for sufficiently large $N$. This shows that, for $N$ sufficiently large, $\mathcal{Q}^N(t)$ will hit the boundary at a time inside $[t - \delta, t)$. Hence, for $N$ sufficiently large,

$$\begin{aligned}
\mathcal{K}(t) &= \lim_{N\to\infty} \mathcal{K}^N(t) \\
&= \lim_{N\to\infty} \left(1 - L\left(\sum_{c=1}^{d} \kappa_c^N \frac{m_c^N(t)}{\mathcal{R}_c^N(t)}\right)^{-1}\right) \\
&= \left(1 - L\left(\sum_{c=1}^{d} \kappa_c \frac{m_c(t)}{\mathcal{R}_c(t)}\right)^{-1}\right).
\end{aligned}$$

In the second case, $\mathcal{K}(t) = p_{\max} = F(q_{\max})$, so

$$\mathcal{K}(t) = \max\left\{p_{\max}, 1 - L\left(\sum_{c=1}^{d} \kappa_c \frac{m_c(t)}{\mathcal{R}_c(t)}\right)^{-1}\right\}.$$

We have therefore established $\mathcal{K}(t) = F(\mathcal{Q}(t))$ if $\mathcal{Q}(t) < q_{\max}$ and

$$\mathcal{K}(t) = \max\left\{p_{\max}, 1 - L\left(\sum_{c=1}^{d} \kappa_c \frac{m_c(t)}{\mathcal{R}_c(t)}\right)^{-1}\right\} \tag{3.7}$$

if $\mathcal{Q}(t) = q_{\max}$ at all times $t$ where $\mathcal{K}(t)$ doesn't jump. Since the derivatives of $m_c(t)$ and $\mathcal{Q}(t)$ are bounded below, there is a time interval to the right of any jump time free of jumps. Since $\mathcal{K}(t)$ is right continuous, it follows that

$$\mathcal{K}(t) = 1 - L\left(\sum_{c=1}^{d} \kappa_c \frac{m_c(t)}{\mathcal{R}_c(t)}\right)^{-1}$$

at jump times.

*Uniform convergence of* $\mathcal{W}_n^N$. We fix some coordinate $n$. Clearly, $\mathcal{W}_n^N(t) \le a(t)$ for all $0 \le t \le T$, so it follows that $\lambda_n^N(t)$ is uniformly bounded by $\overline{\lambda}(t) = a(t)/T_{\min}$ for all $0 \le t \le T$ and all $N \ge n$. Consequently, $\mathcal{N}_n^N(t) \le \overline{N}_n(t)$, where $\overline{N}_n(t) = \int_0^t \int_0^\infty \mathbf{1}_{[0,\overline{\lambda}(s))}(u)\Upsilon_n(du,ds)$. Hence, if for some trajectory $\omega$ of $\Upsilon_n$, $\overline{N}_n(T) \le m$, then $\mathcal{N}_n^N(T) \le m$ for all $n$ and $N$.

For any trajectory $\omega$ of $\Upsilon_n$, we can solve the system

$$\mathcal{W}_n(t) - w_n = \int_0^t \left[\frac{1}{\mathcal{R}_n(s)}\,ds - \frac{\mathcal{W}_n(s^-)}{2}\,d\mathcal{N}_n(s)\right], \tag{3.8}$$

where $\mathcal{W}_n(t) = w_n$, for $t \le 0$, for $n = 1, \dots$ and

$$\mathcal{N}_n(t) = \int_0^t \int_0^\infty \mathbf{1}_{[0,\lambda_n(s))}(u)\Upsilon_n(du,ds)$$

and

$$\lambda_n(s) = \frac{\mathcal{W}_n(s - \mathcal{R}_n(s))}{\mathcal{R}_n(s - \mathcal{R}_n(s))}\mathcal{K}(s - \mathcal{R}_n(s)).$$

Now consider the solution of (3.8) from jump point to jump point of $\mathcal{N}_n$. Let $T_m$, $m = 1, 2, \dots$, be the jumps of $\mathcal{N}_n$ in $[0,T]$ corresponding to jump points $(X_m, Y_m)$, $m = 1, 2, \dots$, of $\Upsilon_n$ which satisfy $T_m = X_m$ and $Y_m \le \lambda_n(X_m)$. The solution of (3.8) for $t \ge T_m$ is the deterministic additive increase of the window size until time $T_{m+1}$ and this has zero chance of hitting $(X_{m+1}, Y_{m+1})$ which is chosen according to a Poisson process. Hence, with probability one, the trace $(t, \lambda_n(t))$ for $0 \le t \le T$ will avoid the points of $\Upsilon_n(du,ds)$ for a given $\omega$, so we can put a band of width $\varepsilon$ around the trace $(t, \lambda_n(t))$. Now as long as $\lambda_n^N(t)$ lies in this band, then the point processes $\mathcal{N}_n$ and $\mathcal{N}_n^N$ are identical. This will be the case if $(\mathcal{R}_n^N(t))^{-1}$ is sufficiently close to $(\mathcal{R}_n(t))^{-1}$ because the solutions $\mathcal{W}_n^N(t)$ and $\mathcal{W}_n(t)$ starting from the same point $w_n$ will rise together inside the band between jumps and at jumps will be cut in half together. We conclude that, with probability one, $\mathcal{W}_n^N(t)$ converges uniformly to $\mathcal{W}_n(t)$ as $N \to \infty$.

*Equation for the limit $\mathcal{W}_n$.* Since $\mathcal{R}_n^N(s)$, $\mathcal{K}^N(s)$ and $\mathcal{W}_n^N(s)$ converge uniformly to $\mathcal{R}_n(s)$, $\mathcal{K}(s-\mathcal{R}_n(s))$ and $\mathcal{W}_n(s)$, respectively, it follows that $\lambda_n^N(t)$ converges uniformly to

$$\lambda_n(t) = \frac{\mathcal{W}_n(s-\mathcal{R}_n(s))}{\mathcal{R}_n(s-\mathcal{R}_n(s))}\mathcal{K}(s-\mathcal{R}_n(s))$$

on $[0,T]$. It therefore follows that both sides of

$$\begin{aligned}\mathcal{W}_n^N(t) - w_n = \int_0^t \frac{1}{\mathcal{R}_n^N(s)}\,ds \\ - \int_0^t \int_{u=0}^{\infty} \frac{\mathcal{W}_n^N(s^-)}{2}\mathbf{1}_{[0,\lambda_n^N(s))}(u)\Upsilon_n(du,ds)\end{aligned}$$

converge uniformly on $D[0,T]$ to (3.8).

3.3. *Equations for the limit $(\mathcal{W},Q)$.*

*Determining $\overline{\mathcal{W}}_c$.* Since $\mathcal{Q}(t)$ is deterministic, the equations in (3.8) are independent. This, in turn, means that

$$E\overline{\mathcal{W}}_c(s) = \overline{\mathcal{W}}_c(s) = \lim_{N\to\infty}\frac{1}{\kappa_c^N N}\sum_{n=1}^{N}\mathcal{W}_n(s)\chi\{n\in K_c\}.$$

Note that this limit is not taken along a subsequence of $N$.

This essentially follows from the law of large numbers. It suffices to consider the $\mathcal{W}_n(s) = f(\mathcal{W}_n(0),\Upsilon_n)$ defined by (3.8). $f$ is defined on $E = [0,\infty)\times D_{R^+}[0,T]$, where $\mathcal{W}_n(0)\in R^+ = [0,\infty)$ and $\Upsilon_n^\omega \in D_{R^+}[0,T]$, the space of cadlag functions. $E$ is a metric space with metric $m = e\oplus d$, where $e$ is the Euclidean metric on $[0,\infty)$ and $d$ is the Skorohod metric on $D_{R^+}[0,T]$.

We have shown above that, for any initial point $w_0$, the set of trajectories $\Upsilon_n^{\omega_0}$ ($\Upsilon_n^{\omega_0}$ is the point process $\Upsilon_n$ evaluated at the sample point $\omega_0$) such that the associated graph $(t,\lambda_n(t))$ does not hit any of the jump points of $\Upsilon_n^{\omega_0}$ has probability one. $f$ is continuous on this set by the arguments used above. If the graph $(t,\lambda_n(t))$ avoids the points of $\Upsilon_n^{\omega_0}$, then according to [4] or [10], for a point $q = (w_1,\Upsilon_n^{\omega_1})\in E$ close to $p = (w_0,\Upsilon_n^{\omega_0})$, we can find a strictly increasing mapping $\theta$ of $[0,T]$ onto itself with $\sup_{[0,T]}|\theta(t)-t|\le\gamma(\theta)$, where $\gamma(\theta)\to 0$ as $q\to p$ such that $\Upsilon_n^{\omega_1}(\theta(t))$ is uniformly close to $\Upsilon_n^{\omega_0}(t)$. Notice that this means the jumps occur at the same times.

The solution to (3.8) for $\mathcal{W}_n$ and $\lambda_n$ for $p$ and for $(w_1,\Upsilon_n^{\omega_1}(\theta(t)))$ are therefore uniformly close. Hence, at any fixed time $s$ where $\Upsilon_n^{\omega_0}$ has no jumps, we will have $\mathcal{W}_n^{\omega_1}(s)$ and $\lambda_n^{\omega_1}(s)$ uniformly close to $\mathcal{W}_n^{\omega_0}(\theta(s))$ and $\lambda_n^{\omega_0}(\theta(s))$. Since $\theta(s)$ is arbitrarily close to $s$ and since the chance of a jump arbitrarily near time $s$ tends to zero, it follows that, for $q$ sufficiently close to $p$, $\mathcal{W}_n^{\omega_1}(s)$ is arbitrarily close to $\mathcal{W}_n^{\omega_0}(s)$. This means $f$ is continuous at $(w_0, \Upsilon_n^{\omega_0})$ as long as $\Upsilon_n^{\omega_0}$ has no jumps at $s$ and this has probability one.

By hypothesis,

$$\lim_{N\to\infty} \frac{1}{\kappa_c^N N} \sum_{n=1}^{N} \delta_{\mathcal{W}_n(0)} \chi\{n \in K_c\} \to \mu_c$$

and the $\Upsilon_n$ are i.i.d. independent Poisson processes, so the empirical measure of the pairs $(\mathcal{W}_n(0), \Upsilon_n)$ converges; that is,

$$\frac{1}{\kappa_c^N N} \sum_{n=1}^{N} \delta_{(\mathcal{W}_n(0),\Upsilon_n)} \chi\{n \in K_c^N\} \to \mu_c \otimes \nu,$$

where $\nu$ is the distribution of $\Upsilon_n$ on $D_{R^+}[0,T]$. The result now follows since $f$ is bounded and the set of discontinuities has probability zero relative to the limiting measure.

*Determining* $\mathcal{M}_c$. In addition, the above means that $\mathcal{M}_c^N(t)$ converges weakly to $\mathcal{M}_c(t)$ almost surely $P$. For any continuous function with compact support, define the limiting measure $\mathcal{M}_c(t)$ by

$$\langle g, \mathcal{M}_c(t)\rangle = \lim_{N\to\infty} \frac{1}{\kappa_c^N N} \sum_{n=1}^{N} g(\mathcal{W}_n^N(t)) \chi\{n \in K_c\}.$$

The deterministic limit exists almost surely by the argument above. This also means that $m_c(t) = \langle Id, M_c(t)\rangle$, so when $\mathcal{Q}(t) = q_{\max}$, $K(t)$ satisfies

$$K(t) = \max\left\{p_{\max}, 1 - L\left(\sum_{c=1}^{d} \kappa_c^N \langle Id, \mathcal{M}_c(s)\rangle \frac{(1-\mathcal{K}(t))}{R_c(t)}\right)^{-1}\right\}.$$

Moreover,

$$\begin{aligned}\overline{\mathcal{W}}_c(t) &:= \lim_{N\to\infty} \frac{1}{\kappa_c^N N} \sum_{n=1}^{N} \mathcal{W}_n^N(t) \chi\{n \in K_c\} \\ &= \langle Id, \mathcal{M}_c(t)\rangle.\end{aligned}$$

*Existence of a strong solution.* Hence, along the subsequence $N$, we obtain an almost sure limit point $(\mathcal{W}, \mathcal{Q}, (\mathcal{M}_1, \dots, \mathcal{M}_d))$ which satisfies the modified system: (3.8) and

$$\begin{aligned}&\mathcal{Q}(t) - \mathcal{Q}(0) \\ &\quad = \int_0^t \Bigg[\sum_{c=1}^{d} \kappa_c \frac{\overline{\mathcal{W}}_c(s)}{\mathcal{R}_c(s)} (1 - \mathcal{K}(s)) - L \\ &\qquad + \left(\sum_{c=1}^{d} \kappa_c \frac{\overline{\mathcal{W}}_c(s)}{\mathcal{R}_c(s)} (1 - F(\mathcal{Q}(s))) - L\right)^{-} \chi\{\mathcal{Q}(s) = 0\}\Bigg]\, ds,\end{aligned} \tag{3.9}$$

where

$$
\begin{aligned}
\overline{\mathcal{W}}_c(t) &:= \lim_{N\to\infty} \frac{1}{\kappa_c^N N} \sum_{n=1}^{N} \mathcal{W}_n(t)\chi\{n \in K_c\} \\
&= \langle Id, \mathcal{M}_c(t)\rangle.
\end{aligned}
\tag{3.10}
$$

We call the above a strong solution associated with the derived subsequence. Moreover, (3.10) is true along all sequences. Hence, the solution to (3.9) and (3.8) is a strong solution of the system in Theorem 2.

*Extension to the timeout and slow-start phases.* If we consider the extended system with timeouts and slow-start, we have to define $N_c^A(t)$, the proportion of the connections from class $c$ in congestion avoidance at time $t$. There will be similar proportions $N_c^U(t)$ in timeout and $N_c^S(t)$ in slow-start. The equation for the queue (neglecting boundary terms) becomes

$$
\frac{dQ^N(t)}{dt} = \sum_{c=1}^{d} [\kappa_c^N N_c^A(t)\langle Id, M_c^N(t)\rangle + \kappa_c^N N_c^S(t)\langle Id, H_c^N(t)\rangle] - L,
$$

where $M_c^N(t)$ is the histogram of the window sizes of connections in congestion avoidance and $H_c^N(t)$ is the histogram of the window sizes of connections in slow start.

We can force the queue to be deterministic by considering the modified system (again neglecting boundary terms):

$$
\frac{dQ^N(t)}{dt} = \sum_{c=1}^{d} [\kappa_c^N E(\mathcal{N}_c^A(t)) E\langle Id, \mathcal{M}_c^N(t)\rangle + \kappa_c^N E(\mathcal{N}_c^S(t)) E\langle Id, \mathcal{H}_c^N(t)\rangle] - L.
$$

The window equations for the modified system are uncoupled as before. We can again pick subsequences so that $\mathcal{Q}^N(t)$ converges and then further subsequences so that $E(N_c^A(t))$, $E(N_c^U(t))$ and $E(N_c^S(t))$ converge. As before, the limiting system is in fact a strong solution to the extended system.

## 4. Uniqueness of strong solutions.

We have constructed a strong solution $(\mathbf{W}, \mathbf{R}, Q, M)$ to (2.13) and (2.12). Our approach is to prove $L^1$ convergence of $(\mathbf{W}^N, Q^N)$ to the strong solution on $[0, T]$:

PROPOSITION 1. *Under Assumptions 1 and 2, $D_N(T) \to 0$ as $N \to \infty$, where*

$$
D_N(t) := E \sup_{\tau \le t} |Q^N(\tau) - Q(\tau)| + E\|\mathbf{W}^N(t) - \mathbf{W}(t)\|,
$$

*where*

$$
\|\mathbf{W}^N(t) - \mathbf{W}(t)\| = \frac{1}{N} \sum_{n=1}^{N} \sup_{0\le\tau\le t} |W_n^N(\tau) - W_n(\tau)|.
$$

We conclude that $Q^N$ and, hence, $R_n^N$ converge in probability to $Q$ and $R_n$, respectively. Hence, for each window $W_n^N$ solving (2.12), the intensity $\lambda_n^N(s)$ converges to

$$\lambda_n(s) = \frac{W_n(s - R_n(s))}{R_n(s - R_n(s))} K(s - R_n(s)).$$

This implies that each window $W_n^N$ converges in probability to $W_n$ in Skorohod norm.

We also have a stronger result than the weak convergence of $\mathbf{W}^N$ to $\mathbf{W}$:

COROLLARY 1. *If Assumptions 1 and 2 hold, then* $\|M_c^N(t) - M_c(t)\|_w \to 0$ *in probability for any* $t \le T$.

PROOF. For any bounded Lipschitz function $g$,

$$\begin{aligned}
&\limsup_{N\to\infty} |E\langle g, M_c^N(t)\rangle - E\langle g, M_c(t)\rangle| \\
&\quad = \limsup_{N\to\infty} \left| E\left[ \frac{1}{\kappa_c^N N} \sum_{n=1}^N g(W_n^N(t))\chi\{n \in K_c\} \right.\right. \\
&\qquad\qquad \left.\left. - \lim_{N\to\infty} \frac{1}{\kappa_c^N N} \sum_{n=1}^N g(W_n(t))\chi\{n \in K_c\} \right] \right| \\
&\quad \le \limsup_{N\to\infty} \left| \frac{1}{\kappa_c^N N} \sum_{n=1}^N E[g(W_n^N(t)) - g(W_n(t))]\chi\{n \in K_c\} \right| \\
&\qquad + \limsup_{N\to\infty} \left| \frac{1}{\kappa_c^N N} \sum_{n=1}^N g(W_n(t))\chi\{n \in K_c\} - E\langle g, M_c(t)\rangle \right| \\
&\quad \le \limsup_{N\to\infty} \frac{1}{\kappa_c^N N} \sum_{n=1}^N E|g(W_n^N(t)) - g(W_n(t))| \\
&\quad \le C_g \limsup_{N\to\infty} \frac{1}{N} \sum_{n=1}^N E|W_n^N(t) - W_n(t)|,
\end{aligned}$$

where $C_g$ is the Lipschitz constant divided by $\min\{\kappa_c^N\}$. Hence,

$$\begin{aligned}
&\limsup_{N\to\infty} |E\langle g, M_c^N(t)\rangle - E\langle g, M_c(t)\rangle| \\
&\qquad \le C_g \limsup_{N\to\infty} \frac{1}{N} \sum_{n=1}^N E \sup_{\tau \le t} |W_n^N(\tau) - W_n(\tau)| \\
&\qquad \le C_g \limsup_{N\to\infty} E\|\mathbf{W}^N(t) - \mathbf{W}(t)\|
\end{aligned}$$

$\to 0.$

Since we can construct a convergence determining sequence based on positive bounded Lipschitz functions, the result follows immediately. $\square$

A similar argument shows that $\|M_c^N(t - R_c^N(t), \cdot; t, \cdot) - M_c(t - R_c^N(t), \cdot; t, \cdot)\|_w \to 0$ in probability for any $t \le T$.

To prove $D_N(t) \to 0$ as $N \to \infty$, we will establish a Gronwall inequality: $D_N(t) \le B_N(t) + C \int_0^t D_N(s)\, ds$, where $B_N(t) \to 0$ as $N \to \infty$. It is easier to establish the Gronwall inequality for Gentle RED where the drop probability function $F_\delta(q)$ rises from $p_{\max}$ at $q_{\max}$ to 1 at $q_{\max} + \delta$ and, hence, is Lipschitz. For Gentle RED, the results in Section 4.1 apply with $\rho = T$. Moreover, we will see below that $B^N(t) = E \int_0^t |\overline{\mathbf{S}}_N(s) - \overline{\mathbf{S}}(s)|\, ds$, where $\overline{\mathbf{S}}_N$ and $\overline{\mathbf{S}}$ were defined by (2.14) and (2.15). This means we can even get a rate of convergence using the Gronwall inequality since $D_N(t) \le B^N(t) + C \int_0^t B^N(s) \exp(C(t-s))\, ds$.

Unfortunately for RED, when $Q^N$ hits the boundary $q_{\max}$, the dynamics change because $F$ is not Lipschitz at $q_{\max}$. Our solution to proving Proposition 1 is to prove convergence on $[0, \rho]$, where $\rho$ to the (deterministic) time when $Q(t)$ first hits $q_{\max}$. This is done in Section 4.1. In Section 4.2 we show the convergence of the transmission rates of the prelimit extends for a time $T_{\min}$ beyond $\rho$ (and, in fact, $T_{\min}$ beyond any point in time). This holds because the transmission rates are determined one RTT in the past. This allows us to extend our proof to cross the boundary when $Q$ hits $q_{\max}$. Then, in Section 4.3, we prove convergence on the interval $[0, \sigma]$, where $\rho < \sigma$ and $\sigma$ is the first time $Q(t)$ leaves the boundary; that is, $Q(t) < q_{\max}$ for $t \in (\sigma, \sigma + \delta]$ for some $\delta > 0$.

We now prove a couple lemmas we will need.

LEMMA 4. *$\dot{Q}(t) > -L + \delta$, where $\delta > 0$ and $\dot{Q}(t) \le a(t)/T_{\min} - L$.*

PROOF. Taking expectations of (2.12) and following the steps of Lemma 2, we get

$$EW_n(t) - w_n \ge \frac{1}{T_{\max} + q_{\max}/L} - \frac{EW_n(t^-)}{2} \frac{a(t)}{T_{\min}}.$$

Since the above bound is positive when $EW_n(t)$ is small, it follows that $EW_n(t)$ is uniformly bounded away from zero for all $t$ and $n$. It therefore follows that $\overline{W}_c(t)$ is uniformly bounded away from zero for all $t$, $n$ and $c$. From (2.13), this means that $\dot{Q}(t) > -L$.

The second inequality follows immediately from (2.13). $\square$

LEMMA 5. *For $0 \le s \le t \le T$, we have that*

$$\sup_{0\le s\le t} |R_n^N(s) - R_n(s)|, \qquad \sup_{0\le s\le t} |R_n^N(s - R_n^N(s)) - R_n(s - R_n(s))|$$

*and*

$$\sup_{0\le s\le t} |Q^N(s - R_n^N(s)) - Q(s - R_n(s))|$$

*are bounded by* $C \sup_{0\le s\le t-T_n} |Q^N(s) - Q(s)|$.

PROOF. Let $\phi_n(s) = T_n + Q(s)/L$, $g_n(s) = s + \phi_n(s)$, $\phi_n^N(s) = T_n + Q^N(s)/L$ and $g_n^N(s) = s + \phi_n^N(s)$. For any $t$, let $g_n(u) = t$ and let $g_n^N(u^N) = t$. Note that $u \le t - T_n$ and $u^N \le t - T_n$. Hence, $R_n^N(t) - R_n(t) = \phi_n^N(u^N) - \phi_n(u) = u^N - u$. Suppose $u \le u^N$ (or $u^N \le u$), the function $g_n$ increases (or decreases) from $g_n(u) = t = g_n^N(u^N)$ to $g_n(u^N)$ [from $g_n(u^N)$ to $g_n^N(u^N)$] by an amount of at least $\delta(u^N - u)/L$ since, by Lemma 4, the derivative of $\phi_n(u)$ is bounded below by $\delta/L$. It follows that $\frac{\delta}{L}|u^N - u| \le |g_n^N(u^N) - g_n(u)| \le \sup_{0\le s\le t-T_n} |g_n^N(s) - g_n(s)| \le \sup_{0\le s\le t-T_n} |Q^N(s) - Q(s)|/L$. Hence,

$$|R_n^N(t) - R_n(t)| = |u^N - u| \le \sup_{0\le s\le t-T_n} |Q^N(s) - Q(s)|/\delta$$

and this gives the first result.

In the same manner we obtained (2.5), we have

$$(1 - \dot{R}_n(t)) = \frac{1}{1 + \dot{Q}(t - R_n(t))/L}$$

and from Lemma 4,

$$\begin{aligned}\dot{R}_n(s) &= \frac{\dot{Q}(s - R_n(s))/L}{1 + \dot{Q}(s - R_n(s))/L} \\ &\le |a(s)/T_{\min} - L|/(\delta/L).\end{aligned}$$

Consequently, using the mean value theorem for $s \le t$, we have

$$\begin{aligned}|R_n(s - R_n^N(s)) - R_n(s - R_n(s))| &\le C|R_n^N(s) - R_n(s)| \\ &\le C \sup_{0\le u\le t-T_n} |Q^N(u) - Q(u)|.\end{aligned}$$

Finally,

$$\begin{aligned}&\sup_{0\le s\le t} |R_n^N(s - R_n^N(s)) - R_n(s - R_n(s))| \\ &\quad \le \sup_{0\le s\le t} |R_n^N(s - R_n^N(s)) - R_n(s - R_n^N(s))|\end{aligned}$$

$$
\begin{aligned}
&\quad + \sup_{0\le s\le t} |R_n(s - R_n^N(s)) - R_n(s - R_n(s))| \\
&\le \sup_{0\le s\le t} |R_n^N(s) - R_n(s)| \\
&\quad + \sup_{0\le s\le t} |R_n(s - R_n^N(s)) - R_n(s - R_n(s))| \\
&\le C \sup_{0\le u\le t-T_n} |Q^N(u) - Q(u)|
\end{aligned}
$$

by the first result and the inequality above.

The third result follows because

$$
\begin{aligned}
&\sup_{0\le s\le t} |Q^N(s - R_n^N(s)) - Q(s - R_n(s))| \\
&\qquad \le \sup_{0\le s\le t} |Q^N(s - R_n^N(s)) - Q(s - R_n^N(s))| \\
&\qquad\quad + \sup_{0\le s\le t} |Q(s - R_n^N(s)) - Q(s - R_n(s))| \\
&\qquad \le \sup_{0\le s\le t-T_n} |Q^N(s) - Q(s)| + C|R_n^N(s) - R_n(s)| \\
&\qquad\qquad\qquad \text{since the derivative of } Q \text{ is bounded and } R_n \ge T_n \\
&\qquad \le C \sup_{0\le u\le t-T_n} |Q^N(u) - Q(u)| \qquad \text{by the first result.} \qquad \square
\end{aligned}
$$

4.1. *Convergence away from the boundary.* We start with $Q^N(0) = q(0)$ in the interior:

- $0 < q(0) < q_{\max}$.

Define $\rho$, respectively $\rho^N$, to be the stopping time when $Q(t)$, respectively $Q^N(t)$, first hits $q_{\max}$. We must define the distance between the marginal process $\mathbf{W}^N(t) \equiv (W_1^N(t), \ldots, W_N^N(t))$, $Q^N(t)$ and $M^N(t) \equiv (M_1^N(t), \ldots, M_d^N(t))$ and the limit processes up to the stopping time $\rho \wedge \rho^N$. For any $t \le \rho$, define

$$
\|\mathbf{W}^N(t) - \mathbf{W}(t)\| = \frac{1}{N} \sum_{n=1}^{N} \sup_{0\le\tau\le t\wedge\rho^N} |W_n^N(\tau) - W_n(\tau)|,
$$

where $\tau$ is a stopping time with respect to $\mathcal{F}_t$. Define

$$
D_N(t) := E \sup_{\tau\le t\wedge\rho^N} |Q^N(\tau) - Q(\tau)| + E\|\mathbf{W}^N(t) - \mathbf{W}(t)\|.
$$

We will establish a Gronwall inequality:

Proposition 2. *$D_N(t) \le B^N(t) + C\int_0^t D_N(s)\,ds$ for $t \le \rho$ and $\sup_{t\le\rho} B^N(t) \to 0$ as $N \to \infty$, where $C$ is a canonical constant throughout this calculation (which unfortunately depends on $F'$). Moreover, $\rho^N$ converges to $\rho$ in probability.*

We just group all universal constants which do not depend on $N$ into $C$. We note that one of the factors in $C$ is the Lipschitz constant, $F'(q), q < q_{\max}$.

*Estimate for $Q$.*

$$\begin{aligned}
&E \sup_{0\le\tau\le t\wedge\rho^N} |Q^N(\tau) - Q(\tau)| \\
&\quad + E \sup_{0\le\tau\le t\wedge\rho^N} \int_0^\tau |\overline{\mathbf{S}}_N^N(s)(1-K^N(s)) - \overline{\mathbf{S}}(s)(1-K(s))|\,ds \\
&\qquad \le E \sup_{0\le\tau\le t\wedge\rho^N} \int_0^\tau |(\overline{\mathbf{S}}_N^N(s) - \overline{\mathbf{S}}_N(s))(1-K^N(s))|\,ds \\
&\qquad\quad + E \sup_{0\le\tau\le t\wedge\rho^N} \int_0^\tau |(\overline{\mathbf{S}}_N(s) - \overline{\mathbf{S}}(s))(1-K^N(s))|\,ds \\
&\qquad\quad + E \sup_{0\le\tau\le t\wedge\rho^N} \int_0^\tau |\overline{\mathbf{S}}(s)(K^N(s) - K(s))|\,ds
\end{aligned} \tag{4.1}$$

and

$$\begin{aligned}
&E \sup_{0\le\tau\le t\wedge\rho^N} |Q^N(\tau) - Q(\tau)| \\
&\quad \le E \sup_{0\le\tau\le t\wedge\rho^N} \int_0^\tau \left( \left| \frac{1}{N} \sum_{n=1}^N \left[ \frac{W_n^N(s)}{R_n^N(s)} - \frac{W_n(s)}{R_n(s)} \right] \right| \right) ds \\
&\qquad + E \sup_{0\le\tau\le t\wedge\rho^N} \int_0^\tau |\overline{\mathbf{S}}_N(s) - \overline{\mathbf{S}}(s)|\,ds \\
&\qquad + E \sup_{0\le\tau\le t\wedge\rho^N} \int_0^\tau \overline{\mathbf{S}}(s)|K^N(s) - K(s)|\,ds \\
&\quad \le \int_0^t \frac{1}{N} \sum_{n=1}^N E\left( \sup_{0\le\tau\le s\wedge\rho^N} \left| \frac{W_n^N(\tau)}{R_n^N(\tau)} - \frac{W_n(\tau)}{R_n(\tau)} \right| \right) ds \\
&\qquad + B^N + E\left( \sup_{0\le\tau\le t\wedge\rho^N} \int_0^\tau \frac{1}{T_{\min}} a(s)|K^N(s) - K(s)| \right) ds,
\end{aligned}$$

where $B^N = E\int_0^t |\overline{\mathbf{S}}_N(s) - \overline{\mathbf{S}}(s)|\,ds$.

Next,

$$
\begin{aligned}
&E\biggl(\sup_{0\le\tau\le s\wedge\rho^N}\biggl|\frac{W_n^N(\tau)}{R_n^N(\tau)} - \frac{W_n(\tau)}{R_n(\tau)}\biggr|\biggr)\\
&\quad\le E\biggl(\sup_{0\le\tau\le s\wedge\rho^N}|W_n^N(\tau) - W_n(\tau)|\frac{1}{R_n^N(\tau)}\biggr)\\
&\qquad + E\biggl(\sup_{0\le\tau\le s\wedge\rho^N}|W_n(\tau)|\biggl|\frac{1}{R_n^N(\tau)} - \frac{1}{R_n(\tau)}\biggr|\biggr)\\
&\quad\le \frac{1}{T_{\min}}E\Bigl(\sup_{0\le\tau\le s\wedge\rho^N}|W_n^N(\tau) - W_n(\tau)|\Bigr)\\
&\qquad + \frac{a(s)}{T_{\min}^2}E\Bigl(\sup_{0\le\tau\le s\wedge\rho^N}|R_n^N(\tau) - R_n(\tau)|\Bigr)\\
&\quad\le \frac{1}{T_{\min}}E\Bigl(\sup_{0\le\tau\le s\wedge\rho^N}|W_n^N(\tau) - W_n(\tau)|\Bigr)\\
&\qquad + \frac{1}{T_{\min}^2}a(s)CE\Bigl(\sup_{0\le\tau\le s\wedge\rho^N}|Q^N(\tau) - Q(\tau)|\Bigr),
\end{aligned}
\tag{4.2}
$$

using Lemma 5.

Moreover,

$$
\begin{aligned}
&E\biggl(\sup_{0\le\tau\le t\wedge\rho^N}\int_0^\tau \frac{1}{T_{\min}}a(s)|K^N(s) - K(s)|\biggr)\,ds\\
&\quad\le CE\biggl(\sup_{0\le\tau\le t\wedge\rho^N}\int_0^\tau |K^N(s) - K(s)|\,ds\biggr)\\
&\quad\le CE\biggl(\sup_{0\le\tau\le t\wedge\rho^N}\int_0^\tau |F(Q^N(s)) - F(Q(s))|\,ds\biggr)\\
&\quad\le C\int_0^t E\Bigl(\sup_{0\le\tau\le s\wedge\rho^N}|Q^N(\tau) - Q(\tau)|\,ds\Bigr).
\end{aligned}
\tag{4.3}
$$

Hence,

$$
\begin{aligned}
&E\Bigl(\sup_{0\le\tau\le t\wedge\rho^N}|Q^N(\tau) - Q(\tau)|\Bigr)\\
&\quad\le \int_0^t \frac{1}{T_{\min}}\frac{1}{N}\sum_{n=1}^N E\Bigl(\sup_{0\le\tau\le s\wedge\rho^N}|W_n^N(\tau) - W_n(\tau)|\Bigr)\,ds
\end{aligned}
$$

$$+\int_0^t \frac{1}{T_{\min}^2} a(s) C E\Big(\sup_{0\le\tau\le s\wedge\rho^N} |Q^N(\tau)-Q(\tau)|\Big)\,ds + B_1^N \tag{4.4}$$

$$+ C\int_0^t E\Big(\sup_{0\le\tau\le s\wedge\rho^N} |Q^N(\tau)-Q(\tau)|\Big)\,ds$$

$$\le B_1^N + C\int_0^t D_N(s)\,ds.$$

*Estimate for $W$.* From (2.12),

$$\frac{1}{N}\sum_{n=1}^N E\Big(\sup_{0\le\tau\le t\wedge\rho^N} |W_n^N(\tau)-W_n(\tau)|\Big)$$

$$\le E\bigg(\sup_{0\le\tau\le t\wedge\rho^N}\int_0^\tau \frac{1}{N}\sum_{n=1}^N E\bigg|\frac{1}{R_n^N(s)}-\frac{1}{R_n(s)}\bigg|\,ds\bigg) \tag{4.5}$$

$$+\frac{1}{2}\frac{1}{N}\sum_{n=1}^N E\bigg(\sup_{0\le\tau\le t\wedge\rho^N}\bigg|\int_0^\tau [W_n^N(s^-)\,dN_n^N(s) - W_n(s^-)\,dN_n(s)]\bigg|\bigg). \tag{4.6}$$

Again, (4.5) is bounded by

$$\int_0^t \frac{a(s)}{T_{\min}^2} E\Big(\sup_{0\le\tau\le s\wedge\rho^N} |R_n^N(\tau)-R_n(\tau)|\Big)\,ds \le C\int_0^t E\Big(\sup_{0\le\tau\le s\wedge\rho^N} |Q^N(\tau)-Q(\tau)|\Big)\,ds. \tag{4.7}$$

By the definition of $N_n^N$,

$$\int_0^\tau \frac{W_n^N(s^-)}{2}\,dN_n^N(s) = \int_{s=0}^\tau\int_{u=0}^\infty \frac{W_n^N(s^-)}{2}\chi_{[0,\lambda_n^N(s))}(u)\Upsilon_n(du,ds).$$

Consequently,

$$\frac{1}{N}\sum_{n=1}^N E\bigg(\sup_{0\le\tau\le t\wedge\rho^N}\bigg|\int_0^\tau [W_n^N(s^-)\,dN_n^N(s) - W_n(s^-)\,dN_n(s)]\bigg|\bigg)$$

$$\le \frac{1}{N}\sum_{n=1}^N E\bigg(\sup_{0\le\tau\le t\wedge\rho^N}\int_0^\tau\int_{u=0}^\infty |W_n^N(s^-)\chi_{[0,\lambda_n^N(s))}(u)$$

$$- W_n(s^-)\chi_{[0,\lambda_n(s))}(u)|\Upsilon_n(du,ds)\Big)$$

$$= \frac{1}{N}\sum_{n=1}^{N} E\Big(\int_0^{t\wedge\rho^N}\int_{u=0}^{\infty} |W_n^N(s^-)\chi_{[0,\lambda_n^N(s))}(u) - W_n(s^-)\chi_{[0,\lambda_n(s))}(u)|\Upsilon_n(du,ds)\Big)$$

$$= \frac{1}{N}\sum_{n=1}^{N} E\Big(\int_0^{t\wedge\rho^N}\int_{u=0}^{\infty} |W_n^N(s^-)\chi_{[0,\lambda_n^N(s))}(u) - W_n(s^-)\chi_{[0,\lambda_n(s))}(u)|\,du\,ds\Big).$$

Hence,

$$\frac{1}{N}\sum_{n=1}^{N} E\Big(\sup_{0\le\tau\le t\wedge\rho^N}\Big|\int_0^{\tau}[W_n^N(s^-)\,dN_n^N(s) - W_n(s^-)\,dN_n(s)]\Big|\Big)$$

$$\le \frac{1}{N}\sum_{n=1}^{N} E\Big(\int_0^{t\wedge\rho^N} |W_n^N(s^-) - W_n(s^-)|\lambda_n^N(s)\wedge\lambda_n(s)\,ds\Big)$$

$$+ \frac{1}{N}\sum_{n=1}^{N} E\Big(\int_0^{t\wedge\rho^N} |W_n^N(s^-)\vee W_n(s^-)|\cdot|\lambda_n^N(s) - \lambda_n(s)|\,ds\Big)$$

$$\le \int_0^t \frac{a(s)}{T_{\min}}\Big[\frac{1}{N}\sum_{n=1}^{N} E\Big(\sup_{0\le\tau\le s\wedge\rho^N}|W_n^N(\tau) - W_n(\tau)|\Big)\Big]\,ds \tag{4.8}$$

$$+ \frac{1}{N}\sum_{n=1}^{N} E\Big(\int_0^{t\wedge\rho^N} a(s)|\lambda_n^N(s) - \lambda_n(s)|\,ds\Big), \tag{4.9}$$

where $\lambda_n^N(s)$ and $\lambda_n(s)$ are less than $(w_n + s/T_{\min})/T_{\min} = a(s)/T_{\min}$.

Also,

$$|\lambda_n^N(s) - \lambda_n(s)|$$

$$\le |W_n^N(s - R_n^N(s)) - W_n(s - R_n^N(s))|\frac{1}{R_n^N(s - R_n^N(s))}K^N(s - R_n^N(s))$$

$$+ |W_n(s - R_n^N(s)) - W_n(s - R_n(s))|\frac{1}{R_n^N(s - R_n^N(s))}K^N(s - R_n^N(s))$$

$$+ |W_n(s - R_n(s))|\frac{1}{R_n^N(s - R_n^N(s))} - \frac{1}{R_n(s - R_n(s))}|K^N(s - R_n^N(s))|$$

$$+ W_n(s - R_n(s))\frac{1}{R_n(s - R_n(s))}|K^N(s - R_n^N(s)) - K(s - R_n(s))|.$$

Hence,

$$|\lambda_n^N(s) - \lambda_n(s)| \le |W_n^N(s - R_n^N(s)) - W_n(s - R_n^N(s))|/T_n \tag{4.10}$$

$$+\,|W_n(s-R_n^N(s))-W_n(s-R_n(s))|/T_n \tag{4.11}$$

$$+\,a(s)|R_n^N(s-R_n^N(s))-R_n(s-R_n(s))|/T_n^2 \tag{4.12}$$

$$+\,a(s)|K^N(s-R_n^N(s))-K(s-R_n(s))|/T_n. \tag{4.13}$$

We must bound $E(\int_0^{t\wedge\rho^N}|\lambda_n^N(s)-\lambda_n(s)|\,ds)$, so we must bound the expectation of the integral of each of the above terms. The first (4.10) satisfies

$$\begin{aligned}E\biggl(\int_0^{t\wedge\rho^N}|W_n^N(s-R_n^N(s))-W_n(s-R_n^N(s))|\,ds/T_n\biggr)\\ \le\frac{1}{T_{\min}}\int_0^t E\biggl(\sup_{0\le\tau\le s\wedge\rho^N}|W_n^N(\tau)-W_n(\tau)|\,ds\biggr).\end{aligned}$$

The second term (4.11) is bounded by

$$\begin{aligned}&E\biggl(\int_0^{t\wedge\rho^N}|W_n(s-R_n^N(s))-W_n(s-R_n(s))|\,ds\biggr)\Big/T_n\\ &\quad\le E\biggl(\int_0^{t\wedge\rho^N}\int_{[s-R_n^N(s)\wedge s-R_n(s),s-R_n^N(s)\vee s-R_n(s)]}\frac{1}{R_n(u)}\,du\,ds\biggr)\\ &\qquad+\frac12E\biggl(\int_0^{t\wedge\rho^N}\int_{[(s-R_n^N(s))\wedge(s-R_n(s)),(s-R_n^N(s))\vee(s-R_n(s))]}\\ &\qquad\qquad W_n(u^-)\,dN_n(u)\,ds\biggr).\end{aligned}$$

Note that

$$\begin{aligned}&\int_0^{t\wedge\rho^N}\chi\{(s-R_n^N(s))\wedge(s-R_n(s))\le u\le(s-R_n^N(s))\vee(s-R_n(s))\}\,ds\\ &\quad=(u+\phi_n^N(u)\vee\phi_n(u))\wedge(t\wedge\rho^N)-(u+\phi_n^N(u)\wedge\phi_n(u))\wedge(t\wedge\rho^N),\end{aligned}$$

where $\phi_n(u)=T_n+Q(u)/L$ and $\phi_n^N(u)=T_n+Q^N(u)/L$.

Hence,

$$\begin{aligned}&E\biggl(\int_0^{t\wedge\rho^N}|W_n(s-R_n^N(s))-W_n(s-R_n(s))|\,ds\biggr)\Big/T_n\\ &\quad\le\frac{1}{T_{\min}}E\biggl(\int_0^{t\wedge\rho^N}|R_n(s)-R_n^N(s)|\,ds\biggr)\\ &\qquad+\frac12E\biggl(\int_0^{t\wedge\rho^N}|\phi_n^N(u)\vee\phi_n(u)-\phi_n^N(u)\wedge\phi_n(u)|W_n(u^-)\lambda_n(u)\,du\biggr)\\ &\quad\le\frac{1}{T_{\min}}E\biggl(\int_0^{t\wedge\rho^N}|R_n(s)-R_n^N(s)|\,ds\biggr)\end{aligned}$$

$$
\begin{aligned}
&+\frac{1}{2}E\biggl(\int_0^{t\wedge\rho^N} |\phi_n^N(u)-\phi_n(u)|\frac{a^2(u)}{T_{\min}}\,du\biggr)\\
&\le \frac{C}{T_{\min}}\int_0^t E\Bigl(\sup_{0\le\tau\le s\wedge\rho^N}|Q(s)-Q^N(s)|\Bigr)\,ds\\
&\quad+\int_0^t \frac{1}{2}\frac{a^2(s)}{T_{\min}L}E\Bigl(\sup_{0\le\tau\le s\wedge\rho^N}|Q^N(\tau)-Q(\tau)|\Bigr)\,ds\\
&\le C\int_0^t E\Bigl(\sup_{0\le\tau\le s\wedge\rho^N}|Q^N(\tau)-Q(\tau)|\Bigr)\,ds,
\end{aligned}
$$

where $C$ is a constant.

To bound the third term (4.12), for $s\le t\wedge\rho^N$,

$$|R_n^N(s-R_n^N(s))-R_n(s-R_n(s))|\le C\sup_{\tau\le s}|Q^N(\tau)-Q(\tau)|$$

by Lemma 5. Taking expectations gives

$$
\begin{aligned}
&E\biggl(\int_0^{t\wedge\rho^N}\frac{a_n(s)}{T_n^2}|R_n^N(s-R_n^N(s))-R_n(s-R_n(s))|\biggr)\,ds\\
&\qquad\le C\int_0^t E\Bigl(\sup_{0\le\tau\le s\wedge\rho^N}|Q^N(\tau)-Q(\tau)|\Bigr)\,ds.
\end{aligned}
$$

Similarly, to bound the fourth term (4.13), for $s\le t\wedge\rho^N$,

$$
\begin{aligned}
&\frac{a(s)}{T_n}|K^N(s-R_n^N(s))-K(s-R_n(s))|\\
&\qquad\le C[|K^N(s-R_n^N(s))-K(s-R_n^N(s))|\\
&\qquad\qquad+|K(s-R_n^N(s))-K(s-R_n(s))|] \qquad (4.14)\\
&\qquad\le C[|Q^N(s-R_n^N(s))-Q(s-R_n^N(s))|\\
&\qquad\qquad+|Q(s-R_n^N(s))-Q(s-R_n(s))|]\\
&\qquad\le C\Bigl[\sup_{\tau\le s}|Q^N(\tau)-Q(\tau)|+|R_n^N(s)-R_n(s)|\Bigr]\\
&\qquad\le C\sup_{\tau\le s}|Q^N(\tau)-Q(\tau)|, \qquad (4.15)
\end{aligned}
$$

since $Q$ has a bounded derivative. Taking expectations shows the fourth term is bounded by

$$C\int_0^t E\Bigl(\sup_{0\le\tau\le s\wedge\rho^N}|Q^N(\tau)-Q(\tau)|\Bigr)\,ds.$$

Hence, we can bound (4.9) by

$$C\int_0^t \Bigg[E\Big(\sup_{0\le\tau\le s\wedge\rho^N}|Q^N(\tau)-Q(\tau)|\Big) + \frac{1}{N}\sum_{n=1}^N E\Big(\sup_{0\le\tau\le s\wedge\rho^N} E|W_n^N(\tau)-W_n(\tau)|\Big)\Bigg]\,ds \le C\int_0^t D_N(s)\,ds.$$

Putting together (4.8), (4.9) and (4.7), we get

$$\frac{1}{N}\sum_{n=1}^N E\Big(\sup_{0\le\tau\le t\wedge\rho^N}|W_n^N(\tau)-W_n(\tau)|\Big) \le C\int_0^t\Bigg[E\Big(\sup_{0\le\tau\le t\wedge\rho^N}|Q^N(\tau)-Q(\tau)|\Big) + \frac{1}{N}\sum_{n=1}^N E\Big(\sup_{0\le\tau\le t\wedge\rho^N}|W_n^N(\tau)-W_n(\tau)|\Big)\Bigg]\,ds.$$

Hence,

$$E\|\mathbf{W}^N(t)-\mathbf{W}(t)\|\le C\int_0^t D_N(s)\,ds.$$

Finally, add in (4.4) and we get our Gronwall inequality: $D_N(t)\le B^N(t)+C\int_0^t D_N(s)\,ds$.

4.2. *When crossing or grazing the boundary.* The construction of the Gronwall inequality in the previous subsection is fairly standard and is sufficient for proving mean field convergence when Gentle RED is used. This subsection, however, resolves the fundamental problem when $Q(t)$ just grazes $q_{\max}$ when RED is used. For example, if we changed the dynamics on the boundary to cause the queue to rapidly grow when $Q$ hit $q_{\max}$, then $Q^N$ would not converge to $Q$ along many sample paths. Those paths where $Q^N$ just missed $q_{\max}$ would drop, while those that hit $q_{\max}$ would rise.

Fortunately, our system allows us to resolve this problem. Using the delay, we show below that $\sup_{0\le\tau\le\rho+T_{\min}}|\overline{S}_N^N(\tau)-\overline{S}(\tau)|\to 0$ as $N\to\infty$. In other words, we can extend the convergence of the transmission rate one RTT into the future beyond the first time to hit the boundary because the evolution of the windows for $[\rho,\rho+T_{\min}]$ is already determined at time $\rho$. This forces $Q^N$ to follow $Q$ for one RTT after hitting the boundary. Hence, if $\overline{S}(t)(1-p_{\max})>L$ for $t\in(\rho,\rho+\delta)$, where $\delta>0$, then we are assured the prelimit $Q^N$ hits the boundary close to where $Q$ hits the boundary and that $\rho^N$ converges to $\rho$. Moreover, we are assured that, for $N$ large enough, there will be a last exit time $\eta^N$, where $\rho^N \le \eta^N < \rho + \delta$ when $Q^N$ jitters off the boundary and that $E|\eta^N - \rho| \to 0$. We can therefore define $\sigma^N$ unambiguously as the infimum of those times after $\rho + \delta$ that $Q^N$ leaves the boundary.

When $Q$ grazes the boundary at time $\rho$, then $\overline{S}(t)(1 - p_{\max}) < L$ for $t \in (\rho, \rho + \delta]$ for some $\delta > 0$. Consequently, $\sigma = \rho$. This case poses a mathematical difficulty because the prelimit $Q^N$ either hits the boundary at a time close to $\rho$ or else avoids the boundary altogether. It is therefore difficult to define $\rho^N$ and $\sigma^N$. We resolve this problem at the end of this subsection by effectively skipping over $\rho$ and saying $Q$ didn't really hit the boundary and the prelimit will at most spend a vanishingly small time on the boundary.

This extension of the convergence of $\overline{S}_N^N(t)$ to $\overline{S}(t)$ for times $t \in [\rho, \rho + T_{\min}]$ is, in fact, valid for any time. The proof doesn't change. In Section 4.3 we use this fact to extend the convergence to $\sigma$ and then by iteration to $T$. (For future reference, it might even be useful to introduce a delay into a system without delay to obtain this property and then prove weak convergence as the delay tends to zero.)

To show convergence one RTT into the future, for $t \le \rho + T_{\min}$, define

$$H(t) = E\left[\frac{1}{N}\sum_{n=1}^{N} \sup_{0 \le \tau \le t} |W_n^N(\tau) - W_n(\tau)|\right].$$

For $t \le T_{\min}$, we can use the same steps as the estimates (4.5) and (4.6) to get

$$H(t) \le E\left[\frac{1}{N}\sum_{n=1}^{N} \sup_{0 \le \tau \le \rho} |W_n^N(\tau) - W_n(\tau)|\right]$$

$$+ \frac{1}{N}\sum_{n=1}^{N}\left(\int_{\rho}^{\rho+T_{\min}} E\left|\frac{1}{R_n^N(s)} - \frac{1}{R_n(s)}\right| ds\right) \tag{4.16}$$

$$+ \frac{1}{2}\frac{1}{N}\sum_{n=1}^{N} E\left(\sup_{\rho \le \tau \le t}\left|\int_{\rho}^{\tau} [W_n^N(s^-)\, dN_n^N(s) - W_n(s^-)\, dN_n(s)]\right|\right). \tag{4.17}$$

Again, (4.16) is bounded by

$$T_{\min}\frac{a(T)}{T_{\min}^2} E\left(\sup_{0 \le \tau \le \rho + T_{\min}} |R_n^N(\tau) - R_n(\tau)|\right)$$

$$\le CE\left(\sup_{0 \le \tau \le \rho} |Q^N(\tau) - Q(\tau)|\right)$$

by Lemma 5. We have already shown convergence up until $\rho$, so this tends to zero as $N \to \infty$.

We estimate (4.17) as we did (4.6) to obtain the following terms corresponding to (4.8) and (4.9):

$$\frac{1}{N}\sum_{n=1}^{N} E\biggl(\sup_{\rho\le\tau\le t}\biggl|\int_\rho^\tau [W_n^N(s^-)\,dN_n^N(s) - W_n(s^-)\,dN_n(s)]\biggr|\biggr)$$

$$\le \int_\rho^t \frac{a(s)}{T_{\min}}\biggl[\frac{1}{N}\sum_{n=1}^{N} E\Bigl(\sup_{\rho\le\tau\le s}|W_n^N(\tau)-W_n(\tau)|\Bigr)\biggr]\,ds \tag{4.18}$$

$$+\frac{1}{N}\sum_{n=1}^{N}\int_\rho^t E(a(s)|\lambda_n^N(s)-\lambda_n(s)|\,ds). \tag{4.19}$$

We can bound (4.19) using the decomposition given by (4.10), (4.11), (4.12) and (4.13). Since each of these terms involve times more than $T_{\min}$ in the past, it is not hard to see that the integral of term (4.10) is bounded by $CE[\frac{1}{N}\sum_{n=1}^N \sup_{0\le t\le\rho}|W_n^N(t)-W_n(t)|]$ and that term (4.11) is bounded by $CE(\sup_{0\le\tau\le\rho}|Q^N(\tau)-Q(\tau)|)$, as is the integral of (4.12) and (4.13). All of these bounds tend to zero as $N\to\infty$.

We conclude that $H(t)\le B^N(t)+C\int_0^t H(s)\,ds$ for $0\le t\le\rho+T_{\min}$, where $B^N(t)$ again denotes a term which goes to zero as $N\to\infty$. Using the Gronwall inequality, we conclude $H(t)\to 0$ as $N\to\infty$. Hence we have convergence of the windows over $[\rho,\rho+T_{\min}]$. Refining the estimate (4.2), we get

$$\begin{aligned} &E\biggl(\sup_{0\le\tau\le\rho+T_{\min}}\biggl|\frac{W_n^N(\tau)}{R_n^N(\tau)}-\frac{W_n(\tau)}{R_n(\tau)}\biggr|\biggr)\\ &\quad\le\frac{1}{T_{\min}}E\Bigl(\sup_{0\le\tau\le\rho+T_{\min}}|W_n^N(\tau)-W_n(\tau)|\Bigr)\\ &\qquad+\frac{1}{T_{\min}^2}a(s)CE\Bigl(\sup_{0\le\tau\le\rho}|Q^N(\tau)-Q(\tau)|\Bigr),\end{aligned}$$

using Lemma 5. Both these estimates tend to zero as $N\to\infty$, so we have shown $\sup_{0\le\tau\le\rho+T_{\min}}|\overline{S}_N^N(\tau)-\overline{S}_N(\tau)|\to 0$ as $N\to\infty$. This completes the argument.

Once we have convergence of the transmission rate until $\rho+T_{\min}$, the case where $Q$ enters the boundary becomes obvious and, clearly, $\rho^N$ converges to $\rho$.

The case where $Q$ grazes the boundary so $\rho=\sigma$ is also resolved because, in probability,

$$\overline{S}_N^N(t)(1-p_{\max})\to\overline{S}(t)(1-p_{\max})<L\qquad\text{for }\rho\le t\le\rho+T_{\min}.$$

Hence, if $Q^N$ enters the boundary at time $\rho^N$, it leaves almost immediately so $\sigma^N-\rho^N\to 0$ in probability. Consequently, we can continue the iteration described in Section 4.3 until the next time, $\rho_1$, $Q$ really hits the boundary and the contribution to the terms

$$E\left(\int_0^{t\wedge\rho_1^N} |K^N(s) - K(s)|\, ds\right)$$

and

$$E\left(\int_0^{t\wedge\rho_1^N} |K^N(s - R_n^N(s)) - K(s - R_n(s))|\, ds\right)$$

by the integral over the interval $[\rho^N, \eta^N]$ is negligible.

The case where $Q$ stays on the boundary and $\overline{S}(t)(1 - p_{\max}) = L$ for $t \in [\rho, \sigma]$ is also theoretically possible. It poses no problem, however, because $Q$ can be considered to be on or off the boundary so the estimates in the previous subsections apply.

4.3. *Mean-field convergence on the boundary.* We now prove convergence on the interval $t \in [0, \sigma]$. Assume $\sigma > \rho$. In Section 4.2 we showed how to handle the case when the queue grazes the boundary. Define $\sigma^N$ to be the infimum over times greater than $\rho + \delta$ that $Q^N$ is less than $q_{\max}$, where $\delta < \sigma - \rho$.

For any $0 \le t \le \sigma$, redefine

$$\|\mathbf{W}^N(t) - \mathbf{W}(t)\| = \frac{1}{N}\sum_{n=1}^{N} \sup_{0\le\tau\le t\wedge\sigma^N} |W_n^N(\tau) - W_n(\tau)|,$$

where $\tau$ is a stopping time with respect to $\mathcal{F}_t$ and $\sigma^N$ is the end of the first sojourn on the boundary by $Q^N$. Redefine

$$D_N(t) := E \sup_{\tau\le t\wedge\sigma^N} |Q^N(\tau) - Q(\tau)| + E\|\mathbf{W}^N(t) - \mathbf{W}(t)\|.$$

Again, we will establish a Gronwall inequality: $D_N(t) \le B_1^N(t) + C\int_0^t D_N(s)\, ds$ for $t \in [0, \sigma]$ where $\sup_{t\le\sigma} B_1^N(t) \to 0$ as $N \to \infty$.

The calculation is almost the same as in Section 4.1. We need only improve the bounds on the terms (4.3) and (4.9) via (4.14). For $s \le t \wedge \sigma^N$, where $t \le \sigma$, there are three possibilities; both $Q^N(s)$ and $Q(s)$ are away from the boundary or both are on the boundary or one is on the boundary and the other isn't. If $Q^N(s) < q_{\max}$ and $Q(s) < q_{\max}$, then

$$|K^N(s) - K(s)| = |F(Q^N(s)) - F(Q(s))| \le C \sup_{\tau\le s\wedge\sigma^N} |Q^N(\tau) - Q(\tau)|.$$

If $Q^N(s) = Q(s) = q_{\max}$, $K^N(s)$ and $K(s)$ are given by

$$\overline{S}_N^N(s)(1 - K^N(s)) = L \quad\text{and}\quad \overline{S}(s)(1 - K(s)) = L,$$

with $K^N(s), K(s) > p_{\max}$. Note this means

$$\overline{S}_N^N(s) \geq L/(1-p_{\max}) \quad \text{and} \quad \overline{S}(s) \geq L/(1-p_{\max}). \tag{4.20}$$

Hence, if $Q^N(s) = Q(s) = q_{\max}$,

$$\begin{aligned}
|K^N(s) - K(s)| &= |L/\overline{S}_N^N(s) - L/\overline{S}(s)| \\
&\leq \frac{(1-p_{\max})^2}{L} |\overline{S}_N^N(s) - \overline{S}(s)| \\
&\leq C(|\overline{S}_N^N(s) - \overline{S}_N(s)| + |\overline{S}_N(s) - \overline{S}(s)|) \\
&\leq C D_N(s) + |\overline{S}_N(s) - \overline{S}(s)|,
\end{aligned}$$

using the same estimate as (4.1).

Hence, to bound expression (4.3), for $t \leq \sigma$,

$$\begin{aligned}
E&\Big(\sup_{0\leq\tau\leq t\wedge\sigma^N} \int_0^\tau |K^N(s) - K(s)|\, ds\Big) \\
&\leq C\int_0^t D_N(s)\, ds + \int_0^t |\overline{S}_N(s) - \overline{S}(s)|\, ds \\
&\qquad + E\Big(\sup_{0\leq\tau\leq t\wedge\sigma^N} \int_0^\tau [\chi\{Q^N(s) < q_{\max} = Q(s)\}]\, ds\Big) \qquad (4.21)\\
&\qquad + E\Big(\sup_{0\leq\tau\leq t\wedge\sigma^N} \int_0^\tau [\chi\{Q(s) < q_{\max} = Q^N(s)\}]\, ds\Big) \\
&\leq C\int_0^t D_N(s)\, ds + B^N(t) + E|\rho^N - \rho| + E|\rho^N - \eta^N|.
\end{aligned}$$

To bound (4.9), we need to improve our bound on (4.14). As before,

$$\begin{aligned}
&|K^N(s - R_n^N(s)) - K(s - R_n(s))| \\
&\quad \leq |K^N(s - R_n^N(s)) - K(s - R_n^N(s))| \qquad (4.22)\\
&\qquad + |K(s - R_n^N(s)) - K(s - R_n(s))|.
\end{aligned}$$

If $Q(s - R_n^N(s)) < q_{\max}$ and $Q(s - R_n(s)) < q_{\max}$, then since $F(q)$ is Lipschitz for $q < q_{\max}$,

$$\begin{aligned}
&|K(s - R_n^N(s)) - K(s - R_n(s))| \\
&\quad = F(Q(s - R_n^N(s)) - F(Q(s - R_n(s)))| \\
&\quad \leq C|Q(s - R_n^N(s)) - Q(s - R_n(s))| \\
&\quad \leq C|R_n^N(s) - R_n(s)| \qquad \text{because } Q \text{ is differentiable} \\
&\quad \leq C \sup_{\tau\leq s\wedge\sigma^N} |Q^N(\tau) - Q(\tau)|.
\end{aligned}$$

However, if $Q(s - R_n^N(s)) = Q(s - R_n(s)) = q_{\max}$,

$$\begin{aligned}
&|K(s - R_n^N(s)) - K(s - R_n(s))| \\
&\quad = |L/\overline{S}(s - R_n^N(s)) - L/\overline{S}(s - R_n(s))| \\
&\quad \le \frac{(1-p_{\max})^2}{L} |\overline{S}(s - R_n^N(s)) - \overline{S}(s - R_n(s))| \\
&\quad \le C|R_n^N(s) - R_n(s)| \\
&\quad \le C \sup_{\tau \le s \wedge \sigma^N} |Q^N(\tau) - Q(\tau)|.
\end{aligned}$$

We used the fact that $\overline{S}(s) = \sum_{c=1}^d \kappa_c \frac{m_c(s)}{R_c(s)}$ is Lipschitz, as was checked in Section 3.2.

We can therefore bound the integral of the second expression in (4.22):

$$\begin{aligned}
&E\biggl(\int_0^{t\wedge\sigma^N} |K(s - R_n^N(s)) - K(s - R_n(s))|\,ds\biggr) \\
&\quad \le C\int_0^t D_N(s)\,ds \\
&\qquad + E\biggl(\sup_{0\le\tau\le t\wedge\sigma^N} \int_0^\tau [\chi\{Q^N(s - R_n^N(s)) < q_{\max} = Q(s - R_n(s))\}]\,ds\biggr) \\
&\qquad + E\biggl(\sup_{0\le\tau\le t\wedge\sigma^N} \int_0^\tau [\chi\{Q(s - R_n(s)) < q_{\max} = Q^N(s - R_n^N(s))\}]\,ds\biggr).
\end{aligned}$$

However, $Q(s - R_n(s)) < q_{\max} = Q^N(s - R_n^N(s))$ implies $s - R_n(s) < \rho$ and $s - R_n^N(s) \ge \rho^N$; that is, when $s \in (\rho^N + \phi_n^N(\rho^N), \rho + \phi_n(\rho))$. Hence,

$$\begin{aligned}
&\int_0^\tau [\chi\{Q(s - R_n(s)) < q_{\max} = Q^N(s - R_n^N(s))\}]\,ds \\
&\quad \le |\rho^N - \rho| + |\phi_n^N(\rho^N) - \phi_n(\rho)| \\
&\quad \le |\rho^N - \rho| + \frac{1}{L}|Q^N(\rho^N) - Q(\rho)| \\
&\quad \le |\rho^N - \rho| + \frac{1}{L}|Q(\rho^N) - Q(\rho)| + C\sup_{\tau\le s\wedge\sigma^N} |Q^N(\tau) - Q(\tau)| \\
&\quad \le C|\rho^N - \rho| + C\sup_{\tau\le s\wedge\sigma^N} |Q^N(\tau) - Q(\tau)|.
\end{aligned}$$

Moreover, $Q^N(s - R_n^N(s)) < q_{\max} = Q(s - R_n(s))$ can occur when $s - R_n(s) \le \rho$ and $s - R_n^N(s) > \rho^N$ or $\rho^N \le s - R_n^N(s) \le \eta^N$; that is, when $s \in$ $(\rho^N + \phi_n^N(\rho^N), \rho + \phi_n(\rho))$ or when $s \in (\rho^N + \phi_n^N(\rho^N), \eta^N + \phi_n^N(\eta^N))$. Hence,

$$\begin{aligned}
&\int_0^\tau [\chi\{Q^N(s - R_n^N(s)) < q_{\max} = Q(s - R_n(s))\}]\, ds \\
&\qquad \le |\rho^N - \rho| + |\rho^N - \eta^N| + |\phi_n^N(\rho^N) - \phi_n(\rho)| + |\phi_n^N(\eta^N) - \phi_n^N(\rho^N)| \\
&\qquad \le 2|\rho^N - \rho| + 2|\rho^N - \eta^N| + \frac{1}{L}|Q^N(\rho^N) - Q(\rho)| \\
&\qquad\quad + \frac{1}{L}|Q^N(\eta^N) - Q^N(\rho^N)| \\
&\qquad \le 2|\rho^N - \rho| + 2|\rho - \eta^N| + \frac{1}{L}|Q(\eta^N) - Q(\rho)| \\
&\qquad\quad + \frac{2}{L}|Q(\rho^N) - Q(\rho)| + C \sup_{\tau \le s \wedge \sigma^N} |Q^N(\tau) - Q(\tau)| \\
&\qquad \le C(|\rho^N - \rho| + |\rho - \eta^N|) + C \sup_{\tau \le s \wedge \sigma^N} |Q^N(\tau) - Q(\tau)|.
\end{aligned}$$

Adding these terms together, we get

$$\begin{aligned}
&E \int_0^{t \wedge \rho^N} |K(s - R_n^N(s)) - K(s - R_n(s))|\, ds \\
&\qquad \le C(|\rho^N - \rho| + |\rho - \eta^N|) + C \int_0^t D_N(s)\, ds.
\end{aligned}$$

We can easily bound the integral of the first expression in (4.22) using the same estimates we made for (4.21):

$$\begin{aligned}
&E\biggl(\sup_{0 \le \tau \le t \wedge \rho^N} \int_0^\tau |K^N(s - R_n^N(s)) - K(s - R_n^N(s))|\, ds\biggr) \\
&\qquad \le C \int_0^t D_N(s)\, ds + B^N(t) + C(E|\rho^N - \rho| + E|\eta^N - \rho|).
\end{aligned}$$

Adding these extra pieces together, we see $D^N(t) \le B_1^N(t) + C \int_0^t D_N(s)\, ds$, where $B_1^N(t) = 2B^N(t) + C(E|\rho^N - \rho| + E|\eta^N - \rho|)$. The first iteration established that $\rho^N$ and $\eta^N$ converge in probability to $\rho$, so $B_1^N(t) \to 0$ as $N \to \infty$. Consequently, $\sup_{t \le \sigma} D^N(t) \to 0$ and we now have convergence on $[0, \sigma]$.

Using the results in Section 4.2, we can show $\sup_{0 \le \tau \le \sigma + T_{\min}} |\overline{S}_N^N(\tau) - \overline{S}(\tau)| \to 0$ as $N \to \infty$. Hence, if $\overline{S}(t)(1 - p_{\max}) < L$ for $t \in (\sigma, \sigma + \delta)$ where $\delta > 0$, then we are assured the prelimit $Q^N$ leaves the boundary close to where $Q$ leaves the boundary and that $\sigma^N$ converges to $\sigma$. We can now iterate to show convergence up through any sequence of entrances and departures from the boundary. It is conceivable that there may be a limit point where a sequence

of entrance and departure times converge to some time $\Theta$. We can use our theory to prove mean field convergence as close as we want to $\Theta$. Then, again using the delay, we can show $\sup_{0\le\tau\le\Theta+T_{\min}}|\overline{S}_N^N(\tau)-\overline{S}(\tau)|\to 0$ as $N\to\infty$. We can therefore establish convergence beyond $\Theta$. There can never be a last time beyond which we cannot establish the mean field convergence.

**5. RED is a weak limit of Gentle RED.** Define the drop probability function for Gentle RED by $F_\delta(q)$, which rises from $p_{\max}$ at $q_{\max}$ to 1 at $q_{\max}+\delta$. Section 4.1 gives mean field convergence for Gentle RED over any interval $[0,T]$ since $F_\delta$ is Lipschitz. Here we show that, as $\delta\to 0$, Gentle RED becomes RED and $F_\delta$ tends (weakly) to $F$, the drop probability function for RED.

For any $\delta$ and any $N$, redefine the solution to the $N$-particle system in Section 2, $(\mathbf{W}^N(t),Q^N(t))$ by $(\mathbf{W}^{\delta,N},Q^{\delta,N})$, so now $(\mathbf{W}^N(t),Q^N(t))$ only denotes the solution with loss function $F$. These processes are constructed iteratively on the almost surely finite number of segments defined by jumps of $\Upsilon_n(s,\overline{\lambda}(T));n=1,\dots,N$, where $\Upsilon_n,n=1,2,\dots$ are defined on the probability space $(\Omega,\mathcal{F},P)$. Let $\mathbf{R}^{\delta,N}=(R_1^{\delta,N},\dots,R_N^{\delta,N})$ be the corresponding round trip time delay of the connections. Let

$$\overline{S}_N^{\delta,N}(t)=\sum_{c=1}^{d}\kappa_c^N\frac{\overline{\mathbf{W}}_c^{\delta,N}(t)}{R_c^{\delta,N}(s)}.$$

Let $P^{\delta,N}$ be the measure induced on $D[0,T]\times C[0,T]$ by $(\overline{S}_N^{\delta,N},Q^{\delta,N})$ (where coordinates greater than $N$ are identically zero). In the same manner as Lemma 1, we can show the measures $P^{\delta,N}$ are tight.

Using this lemma we can now prove the following:

LEMMA 6. *$(\mathbf{W}^{\delta,N},Q^{\delta,N},\mathbf{R}^{\delta,N})$ converge weakly to $(\mathbf{W}^N(t),Q^N,\mathbf{R}^N)$ as $\delta\to 0$. The lemma holds even for $N=\infty$.*

PROOF. By hypothesis, $Q^{\delta,N}(t)=Q^N(t)=q(0)$ and $W_n^{\delta,N}(t)=W_n^N(t)=w_n$ for $t\le 0$. We show the drop probability $K^{\delta,N}(t)=F_\delta(Q^{\delta,N}(t))$ converges as $\delta\to 0$.

First, pick a subsequence $\delta_k$ such that $P^{\delta_k,N}$ converges weakly to $P^{0,N}$ and, moreover, such that $(\mathbf{W}^{\delta_k,N},\overline{S}_N^{\delta_k,N},Q^{\delta_k,N})$ converges almost surely to $(\mathbf{W}^{0,N},\overline{S}_N^{0,N},Q^{0,N})$.

If $Q^{\delta_k,N}(t)<q_{\max}$, then $K^{\delta_k,N}(t)=F_{\delta_k}(Q^{\delta_k,N}(t))=F(Q^{\delta_k,N}(t))$. On the other hand, if $Q^{\delta_k,N}(t)>q_{\max}$, then there will exist a time $\rho^{\delta_k,N}(t)\le t$ when $Q^{\delta_k,N}(t)$ last hit $q_{\max}$. Consider the solution to (2.9) over the interval $[\rho^{\delta_k,N}(t),t]$ and let $V^{\delta_k,N}(s)=1-F_{\delta_k}(Q^{\delta_k,N}(s))$. Note that, for $\rho^{\delta_k,N}(t)\le$

$s \le t$,

$$\frac{dV^{\delta_k,N}(s)}{ds} = -c_{\delta_k}\overline{S}_N^{\delta_k,N}(s)V^{\delta_k,N}(s) + Lc_{\delta_k},$$

where $c_{\delta_k} = (1-p_{\max})/\delta_k$. Solve this equation from time $\rho^{\delta_k,N}(t)$, where $V^{\delta_k,N}(\rho^{\delta_k,N}(t)) = (1-p_{\max})$ up to time $t$:

$$\begin{aligned} V^{\delta_k,N}(t) &= (1-p_{\max})\exp\biggl(-\int_{\rho^{\delta_k,N}(t)}^t c_{\delta_k}\overline{S}_N^{\delta_k,N}(u)\,du\biggr) \\ &\quad + \int_{\rho^{\delta_k,N}(t)}^t Lc_{\delta_k}\exp\biggl(-\int_u^t c_{\delta_k}\overline{S}_N^{\delta_k,N}(s)\,ds\biggr)du \\ &= (1-p_{\max})\exp(-v(\rho^{\delta_k,N}(t))) + \int_0^{v(\rho^{\delta_k,N}(t))} e^{-v}\frac{L}{\overline{S}_N^{\delta_k,N}(u(v))}\,dv, \end{aligned}$$

where $v = \int_u^t c_{\delta_k}\overline{S}_N^{\delta_k,N}(s)\,ds$ and $u(v)$ is the inverse defined implicitly.

Define $D^{\delta_k,N}(t) = 0$ if $Q^{\delta_k,N}(t) \le q_{\max}$ and for $Q^{\delta_k,N}(t) > q_{\max}$, define

$$\begin{aligned} D^{\delta_k,N}(t) &= V^{\delta_k,N}(t) - \frac{L}{\overline{S}_N^{\delta_k,N}(t)} \\ &= (1-p_{\max})\exp(-v(\rho^{\delta_k,N}(t))) \\ &\quad + L\int_0^{v(\rho^{\delta_k,N}(t))} e^{-v}\biggl(\frac{1}{\overline{S}_N^{\delta_k,N}(u(v))} - \frac{1}{\overline{S}_N^{\delta_k,N}(t)}\biggr)dv \\ &\quad - e^{-v(\rho^{\delta_k,N}(t))}\frac{L}{\overline{S}_N^{\delta_k,N}(t)}. \end{aligned}$$

Since $(\overline{S}_N^{\delta_k,N}(t), Q^{\delta_k,N(t)})$ converges almost surely to $(\overline{S}_N^{0,N}(t), Q^{0,N}(t))$ as $\delta_k \to 0$, it follows that $\rho^{\delta_k,N}(t)$ converges to a limit $\rho^{0,N}(t)$. Since $\overline{S}_N^{0,N}(t) > 0$ almost surely, $v(\rho^{\delta_k,N}(t)) \to \infty$ as $\delta_k \to 0$ and for a fixed $v$, $v\frac{\delta_k}{(1-p_{\max})} = \int_{u(v)}^t \overline{S}_N^{\delta_k,N}(s)\,ds$ so $u(v) \to 0$ as $\delta_k \to 0$. We conclude $D^{\delta_k,N}(t) \to 0$ almost surely.

If we now take the limit as $\delta_k \to 0$ in (2.9) satisfied by $(\overline{S}_N^{\delta_k,N}(t), Q^{\delta_k,N}(t))$, using loss function $F_{\delta_k}$, we see $(\overline{S}_N^{0,N}(t), Q^{0,N(t)})$ satisfies (2.9) with loss function $F$, so $K^{0,N}(t) = L/\overline{S}_N^{0,N}(t)$ if $Q^{0,N}(t) = q_{\max}$. Moreover, the window $W_n^{\delta_k,N}$ satisfies (2.3), where the rate of window reductions

$$\lambda_n^{\delta_k,N}(t) := \frac{W_n^{\delta_k,N}(t - R_n^{\delta_k,N}(t))}{R_n^{\delta_k,N}(t - R_n^{\delta_k,N}(t))}F_{\delta_k}(Q^{\delta_k,N}(t - R_n^{\delta_k,N}(t))).$$

Taking the limit as $\delta_k \to 0$, we see $W_n^{\delta_k,N}$ converges to $W_n^{0,N}$, satisfying (2.3) where the rate of window reductions is

$$\lambda_n^{0,N}(t) := \frac{W_n^{0,N}(t - R_n^{0,N}(t))}{R_n^{0,N}(t - R_n^{0,N}(t))} F_{\delta_k}(Q^{0,N}(t - R_n^{0,N}(t))).$$

All this means $\mathbf{W}^{0,N} = (W_1^{0,N}, \ldots, W_N^{0,N})$ and $Q^{0,N}$ are a solution to (2.3) and (2.9), so, in fact, must be $\mathbf{W}^N = (W_1^N, \ldots, W_N^N)$ and $Q^N$. Hence, the limit along subsequences is unique so we have proved weak convergence. □

**6. Mean-field stochastic differential equations.** In this section we prove Theorem 1. We can reformulate (2.8) as in [2]. For $g \in \mathcal{G}$,

$$\begin{aligned}
&\langle g, M_c^N(t)\rangle - \langle g, M_c^N(0)\rangle \\
&\quad = \frac{1}{\kappa_c^N N} \sum_{n=1}^N \int_0^t \biggl[\frac{dg}{dw}(W_n^N(s)) \frac{1}{R_c^N(s)}\, ds \\
&\qquad\qquad + (g(W_n^N(s^-)/2) - g(W_n^N(s^-)))\, dN_n(s)\biggr] \chi\{n \in K_c\}
\end{aligned}$$

$$\begin{aligned}
&\quad = \frac{1}{\kappa_c^N N} \sum_{n=1}^N \chi\{n \in K_c\} \int_0^t \biggl[\frac{dg}{dw}(W_n(s)) \frac{1}{R_c^N(s)}\, ds \\
&\qquad\qquad + (g(W_n^N(s)/2) - g(W_n^N(s))) \\
&\qquad\qquad\quad \times \frac{W_n^N(s - R_c^N(s))}{R_c^N(s - R_c^N(s))} K^N(s - R_c^N(s))\, ds\biggr] \\
&\qquad + \mathcal{E}_c^N(t),
\end{aligned} \tag{6.1}$$

where $\mathcal{E}_c^N(t)$ is given by

$$\frac{1}{\kappa_c^N N} \sum_{n=1}^N \chi\{n \in K_c\} \int_0^t (g(W_n^N(s^-)/2) - g(W_n(s^-)))\, dZ_n^N(s)$$

and

$$Z_n^N(t) - Z_n^N(0) := \int_0^t \biggl(dN_n^N(s) - \frac{W_n^N(s - R_n^N(s))}{R_c^N(s - R_c^N(s))} K^N(s - R_c^N(s))\, ds\biggr).$$

Hence,

$$\begin{aligned}
&\langle g, M_c^N(t)\rangle - \langle g, M_c^N(0)\rangle \\
&\quad = \int_0^t \biggl[\frac{1}{R_c^N(s)} \biggl\langle \frac{dg(w)}{dw}, M_c^N(s)\biggr\rangle\, ds
\end{aligned}$$

$$+\langle (g(w/2)-g(w))v, M_c^N(s-R_c^N(s),dv;s,dw)\rangle \tag{6.2}$$

$$\times \frac{1}{R_c^N(s-R_c^N(s))}K^N(s-R_c^N(s))\,ds\Bigg]$$

$$+\mathcal{E}_c^N(t).$$

We first show $\mathcal{E}_c^N$ is asymptotically small as $N\to\infty$. Recall

$$\mathcal{E}_c^N(t)=\frac{1}{\kappa_c^N N}\sum_{n=1}^N\int_0^t C_{c,n}^N(s)Z_{c,n}^N(ds),$$

where

$$C_{c,n}^N(s)=\chi\{n\in K_c\}(g(W_n^N(s)/2)-g(W_n^N(s)))$$

and

$$Z_{c,n}^N(t)-Z_{c,n}^N(0) := \int_0^t\left(dN_n^N(s)-\frac{W_n^N(s-R_c^N(s))}{R_c^N(s-R_c^N(s))}K^N(s-R_c^N(s))\,ds\right)\chi\{n\in K_c\}.$$

If $n\in K_c$, $N_n^N(s)$ is a point process adapted to $\mathcal{F}_n(t)$ with a stochastic intensity $W_n^N(s-R_c^N(s))K^N(s-R_c^N(s))/R_c^N(s-R_c^N(s))$. Consequently, $Z_n^N(t)$ is a martingale. Recall that $W_n^N(s)$ is also adapted to $\mathcal{F}_n(t)$, so the right continuous version is $\mathcal{F}_n(t)$-predictable. By Theorem T13 in [1],

$$\begin{aligned}
&E(\mathcal{E}_N^c(t))^2\\
&\quad=\frac{1}{(\kappa_N^c N)^2}\\
&\qquad\times E\left[\sum_{n=1}^N\chi\{n\in K_c\}\int_0^t (C_{n,c}^N)^2(s)\frac{W_n^N(s-R_c^N(s))}{R_c^N(s-R_c^N(s))}K^N(s-R_c^N(s))\,ds\right]\\
&\quad\le\frac{1}{(\kappa_N^c N)^2}\\
&\qquad\times E\sum_{n=1}^N\chi\{n\in K_c\}\int_0^t\left(g(W_n^N(s)/2)-g(W_n^N(s))^2\frac{W_n^N(s-R_c^N(s))}{R_c^N(s-R_c^N(s))}\,ds\right)\\
&\quad\le\frac{C_1}{(\kappa_N^c N)^2}\left(\sum_{n=1}^N\chi\{n\in K_c\}\int_0^t E[W_n^N(s-R_c^N(s))]\,ds\right),
\end{aligned}$$

where $C_1$ is a constant depending on $\sup g$ and $\sup(g)'$.

We have the a priori bound $W_n^N(t) \le a(t)$. Hence,

$$E(\mathcal{E}_c^N(t))^2 \le \frac{tC_1}{(\kappa_c^N)^2 N} a(t).$$

So $\mathcal{E}_c^N(t)$ tends to 0 in $L^2$. Since, in addition, $\mathcal{E}_c^N(t)$ is a martingale, it follows that

$$P\Big(\sup_{t\in[0,T]} |\mathcal{E}_c^N(t)| > \lambda\Big) \le E(\mathcal{E}_c^N(T))^2/\lambda^2,$$

so the process $\mathcal{E}_c^N(t), t \in [0,T]$, converges to zero in probability.

The processes $Q^N(t)$ and $K^N(t)$ converge in probability to the limit processes $Q(t)$ and $K(t)$, while $(M_1^N(t),\dots,M_d^N(t))$ converges to $(M_1(t),\dots,M_d(t))$ in probability where the limit processes satisfy (2.13) and (2.12). Take the limit of (6.2) and we have our proof.

**7. Numerical analysis.** Assuming $g(0)=0$ and that $g(w)\to 0$ as $w\to\infty$, we can rewrite (1.2) as

$$\begin{aligned}
&\langle g, M_c(t)\rangle - \langle g, M_c(0)\rangle \\
&\quad = \int_0^t \bigg[ -\frac{1}{R_c(s)} \langle g(w), D_w M_c(s,dw)\rangle \\
&\qquad\quad + \frac{1}{R_c(s-R_c(s))} K(s-R_c(s)) \\
&\qquad\quad \times \langle g(w), e(s,s-R_c(s),2w)\cdot M_c(s,2dw) \\
&\qquad\qquad - e(s,s-R_c(s),w)\cdot M_c(s,dw)\rangle \bigg]\, ds,
\end{aligned}$$

where $D_w M_c(s,dw)$, respectively $D_t M_c(s,dw)$, is the Frechet derivative of the measure $M_c(s,dw)$ with respect to $w$, respectively $t$. Consequently,

$$\begin{aligned}
D_t M_c(t,dw) &= -\frac{1}{R_c(t)} D_w M_c(t,dw) \\
&\quad + \frac{1}{R_c(t-R_c(t))} K(t-R_c(t)) \\
&\quad \times (e(t,t-R_c(t),2w) M_c(t,d(2w)) \\
&\qquad - e(t,t-R_c(t),w) M_c(t,dw)).
\end{aligned} \tag{7.1}$$

Neither $M_c(t,dw)$ or $M_c(s-R_c(s),dv;s,dw)$ is a state, but the above equation does provide enough information to evolve the system. Let $\mu_c(t)$ denote the process $\{M_c(s, dw); t-1 \le s \le t\}$ (all RTTs are less than 1). Using (1.1), we can evolve $M_c(t, dw)$ from $t$ to $t+\delta t$ while $M_c(t-s+\delta t, dw)$ is obtained by a time shift. Unfortunately, $\mu_c$ is not a practical state. Even if we discretize and only keep the trajectory of the process on a partition giving $\{M_c(s_i, dw); t-1 = s_0 < s_1 < \cdots < s_n = t\}$, it still requires too much computer memory to solve numerically.

We can avoid this problem by defining a sequence of times $t_k^c$ for each class such that $t_{k+n+1}^c - R_c(t_{k+n+1}^c) = t_k^c$. If we pick $n$ sufficiently large, this gives a fine partition. Define $\Phi^c(t) =$ the first $k$ such that $t_k^c > t$. We will construct our solution by recurrence from time $t_i$ to $t_{i+1}$ by defining $t_{i+1} = \min_c(t_{\Phi^c(t_i)}^c)$ starting from time $t_0 = 0$. Next assume that, for each class, we have been able to calculate and save the vector $V_c^M(t)$, a discretized version of $M_c(t_k^c)$ for $k = m-n, \ldots, m$, where $m = \Phi^c(t) - 1$ (these are marginals, not the entire joint distribution). Also assume we save the vector of kernels $V_c^T(t)$ given by $T^c(t_k^c)$ for $k = m-n, \ldots, m$, where $M_c(t_k^c) = M_c(t_{k-1}^c) \circ T^c(t_k^c)$ and $m$ is as above. Finally, assume that we save the kernels $S^c(m) = \prod_{k=m-n}^m T^c(t_k^c)$.

We can now evolve our system to $t_{i+1}$. At each step, we evolve the queue and the one class $c$, where $t_{i+1} = t_{\Phi^c(t_i)}^c$. The inverse kernel $(S^c(m))^{-1}$, gives the conditional distribution of the windows of class $c$ one RTT before time $t_i$, given the window at time $t_i$. Calculate the conditional expectation $e(t_i, t_i - R_c(t_i), w)$. With this we can use (1.2) to calculate $T^c(t_{m+1}^c)$. Drop $T^c(t_{m-n}^c)$. Update $S^c(m+1) = (T^c(t_{m-n}^c))^{-1} S^c(m) T^c(t_{m+1}^c)$. Finally, we calculate $Q(t_{i+1})$ using the $M_c(\Phi^c(t_i) - 1)$ for $c = 1, \ldots, d$.

In Section 3.3.2 in [2] we made a smooth approximation to $e(t, t-R_c(t), w)$ based on the fact that one RTT in the past the window size was most likely the current window size minus once or twice the current window size if a loss was detected in the interim. With this approximation, we used (7.1) to evolve a discrete approximation of the measure $M_c$ (except [2] only treats one class). The numerical results are excellent after one corrects for the fact that a proportion of the connections in an Opnet simulation are in timeout (our model assumes connections instantaneously resume congestion avoidance if they fall into timeout).

To illustrate the mean field limit, we performed an Opnet simulation with $N = 200$, $N = 400$ and $N = 800$ sources (see Figures 1–4). Each source sends packets of size 536 bytes to a T3 router with a transmission rate of 44.736 Megabits per second or $L = 10433$ packets per second. We assume the sources all have a transmission delay of 100 milliseconds. The router implements RED with $p_{\max} = 0.05$ for all the simulations, but we rescale $Q_{\max}$ to be 1000 with 200 sources, 2000 with 400 sources and 4000 with 800 sources. Since $Q_{\max}$ scales with $N$, the average queue size does as well while holding the loss probability fixed, which in turn holds the average window size fixed. As $N$ increases we see the fluctuations in the relative queue size (in packets per connection) decrease. We also see that the relative queue size of the Matlab numerical simulation is a bit high. This is because of timeouts as discussed in [2].

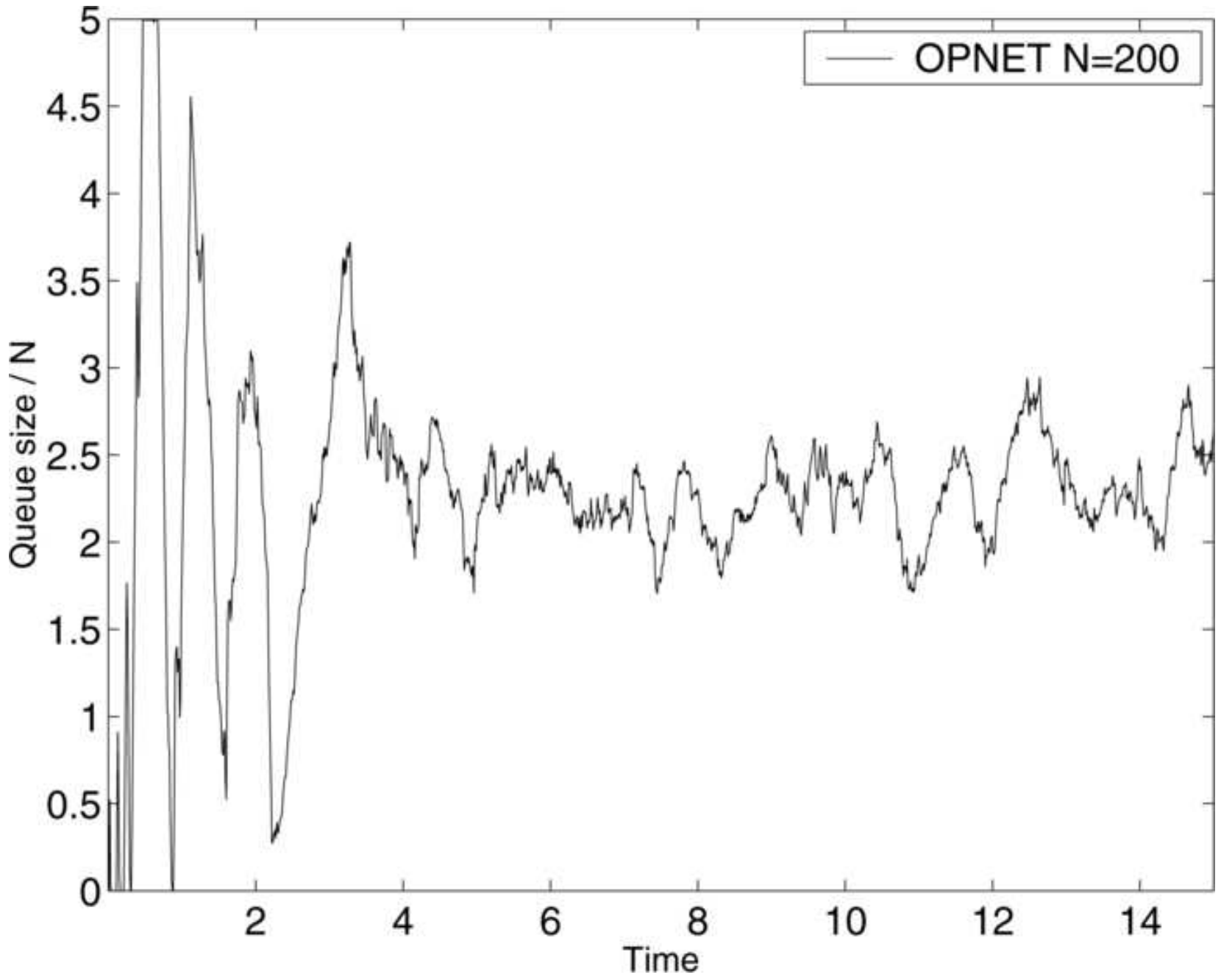

FIG. 1. *Relative queue size with* 200 *sources.*

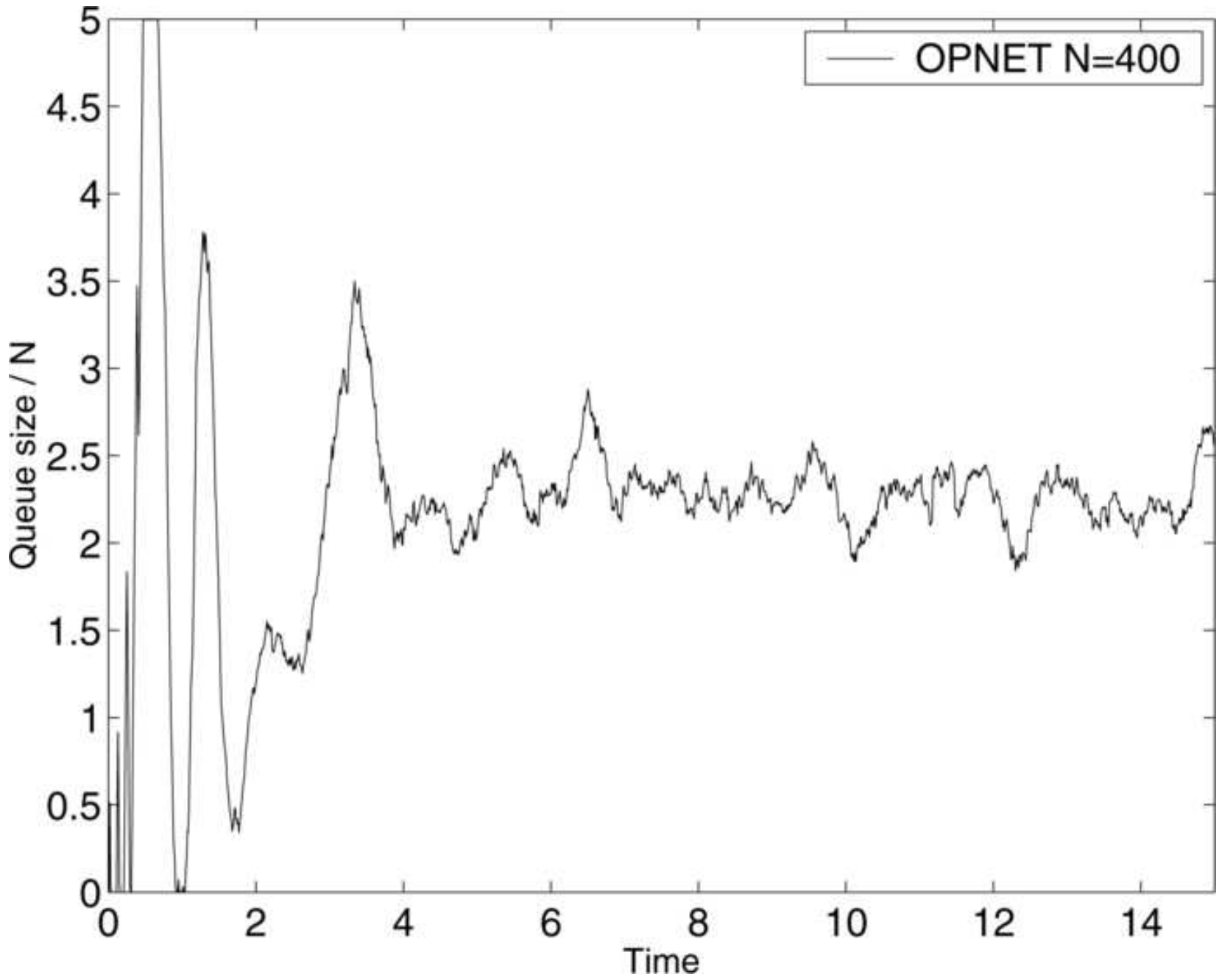

FIG. 2. *Relative queue size with* 400 *sources.*

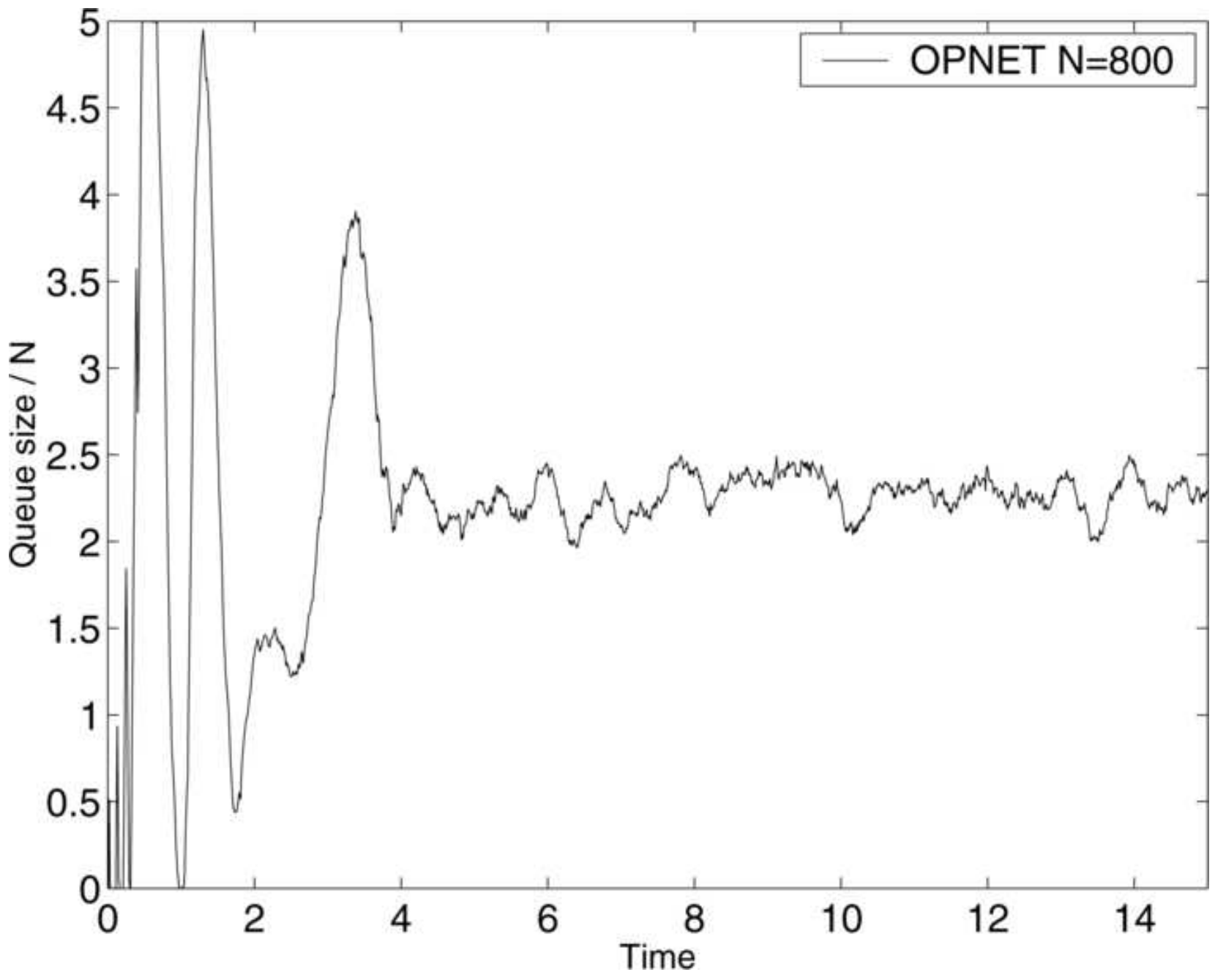

Fig. 3. *Relative queue size with* 800 *sources.*

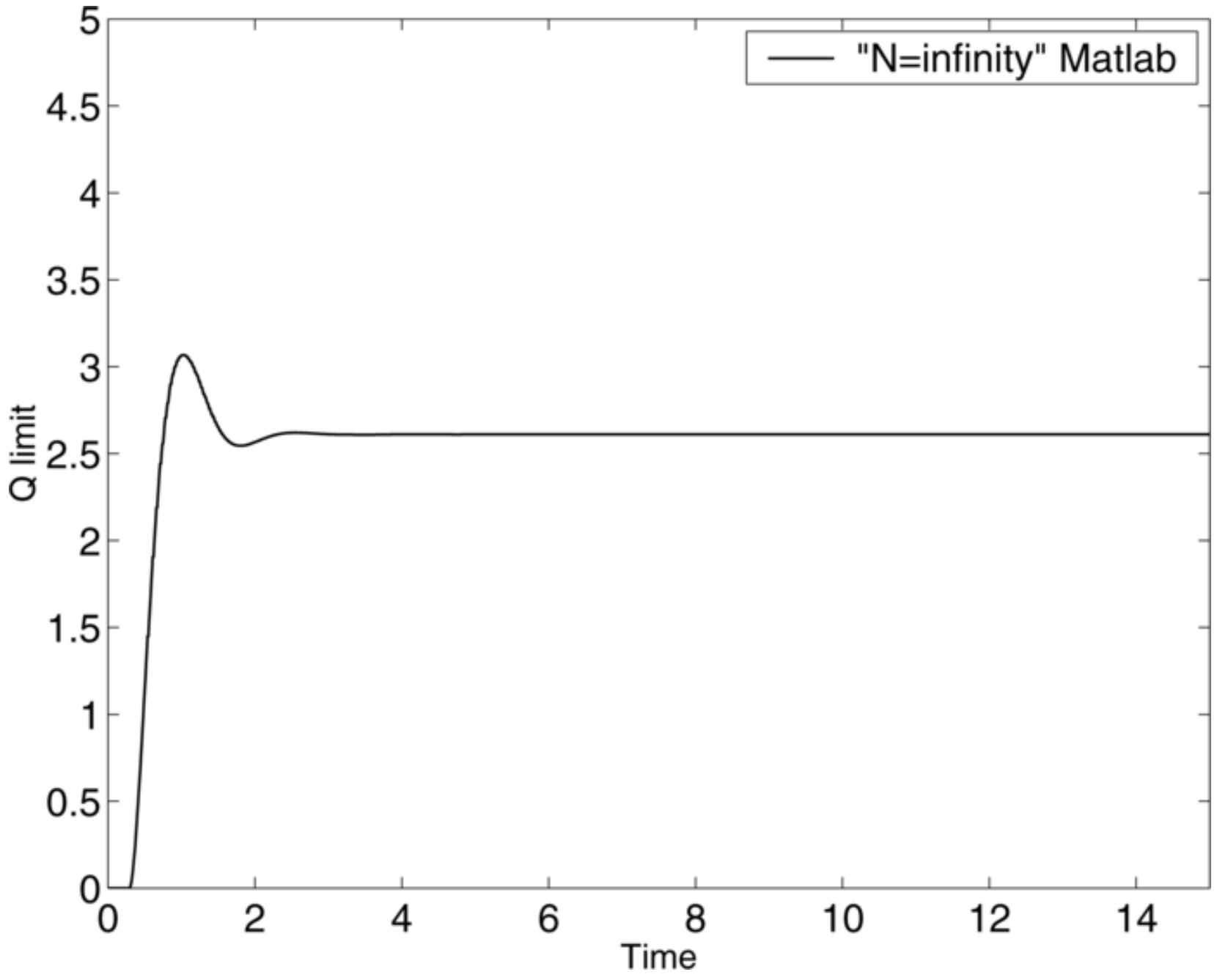

Fig. 4. *Queue size of the mean field limit.*

**Acknowledgments.** J. Reynier thanks François Baccelli for the supervision of this work. R. McDonald thanks Tom Kurtz for his tutorial on how to *bring back the particles* into the proof of weak convergence to the mean-field model. He also thanks Ruth Williams for her suggestions [14]. Both of us thank Pierre Brémaud for stopping us from going completely off the rails and finally we thank a careful referee for his many comments and suggestions.

Department of Mathematics
University of Ottawa
Canada
e-mail: dmdsg@mathstat.uottawa.ca

Département d'Informatique
École Normale Supérieure
France
e-mail: Julien.Reynier@ens.fr